\documentstyle[12pt]{article}

\input{emlines.sty}

\newcommand{\mn}{\sl}
\def\afterthmseparator{}
\makeatletter
\renewcommand{\@begintheorem}[2]{\trivlist
      \item[\hskip \labelsep{\bf #1\ #2\unskip\afterthmseparator}]\mn}
\renewcommand{\@opargbegintheorem}[3]{\trivlist
      \item[\hskip \labelsep{\bf #1\ #2\ (#3)\unskip\afterthmseparator}]\mn}
\makeatother
\newtheorem{theorem}{Theorem}[section]
\newtheorem{lemma}[theorem]{Lemma}
\newtheorem{corollary}[theorem]{Corollary}
\newtheorem{proposition}[theorem]{Proposition}
\newtheorem{rem}[theorem]{Remark}
\newenvironment{remark}{\renewcommand{\mn}{\rm} \begin{rem}}{\end{rem}}
\newtheorem{probl}[theorem]{Problem}

\newtheorem{df}[theorem]{Definition}
\newenvironment{definition}{\renewcommand{\mn}{\rm} \begin{df}}{\end{df}}
\newtheorem{exmpl}[theorem]{Example}

\newcounter{oq}
\newcommand{\que}{\refstepcounter{oq}\par{\bf \theoq.}~}

\newcommand{\bull}{\mbox{$\;\;\;$\vrule height .9ex width .8ex depth -.1ex}}
\newcommand{\qed}{\mbox{$\;\;\;\Box$}}
\newcommand{\noproof}{\unskip\nobreak\hfill \bull}
\newenvironment{proof}{\par\smallbreak\noindent{\sl Proof.~}}%
{\unskip\nobreak\hfill \bull \par\medbreak}
\newenvironment{proofof}[1]{\par\smallbreak\noindent{\sl Proof~of~#1.~}}%
{\unskip\nobreak\hfill \bull \par\medbreak}
{\qed \par\smallbreak}
\newcounter{claim}[theorem]
\renewcommand{\theclaim}{\thetheorem.\arabic{claim}}
{\par\smallskip\par}
\newcommand{\case}[2]{\par{\it Case #1:\/ #2.}}
\newcommand{\subcase}[2]{\par{\it Subcase #1:\/ #2.}}

\newcommand{\hide}[1]{}
\newcommand{\refeq}[1]{(\ref{#1})}
\newcommand{\function}[2]{:#1 \rightarrow #2}
\newcommand{\of}[1]{\left( #1 \right)}
\newcommand{\ofcurle}[1]{\left\{ #1 \right\}}
\newcommand{\setdef}[2]{\left\{\hspace{0.5mm}#1:\hspace{0.5mm} #2\right\}}

\newcommand{\ga}{\alpha}
\newcommand{\gb}{\beta}
\newcommand{\gs}{\sigma}
\newcommand{\And}{\wedge}
\newcommand{\Or}{\vee}

\newcommand{\ra}{\rightarrow}

\newcommand{\ah}{\alpha}

\newcommand{\ben}{\begin{enumerate}}
\newcommand{\een}{\end{enumerate}}
\newcommand{\bec}{\begin{center}}
\newcommand{\enc}{\end{center}}

\newcommand{\beq}{\begin{equation}}
\newcommand{\eeq}{\end{equation}}

\newcommand{\EF}{Ehrenfeucht}
\newcommand{\BS}{Bernays-Sch\"onfinkel}
\newcommand{\game}{\mbox{\sc Ehr}}
\newcommand{\Nest}{\mathop{\sl Nest}\nolimits}
\newcommand{\nest}[1]{\Nest(#1)}
\newcommand{\Qr}{\mathop{\sl qr}\nolimits}
\newcommand{\qr}[1]{\Qr(#1)}
\newcommand{\Alt}{\mathop{\sl alt}\nolimits}
\newcommand{\alt}[1]{\Alt(#1)}
\newcommand{\altset}[1]{\Lambda_{#1}}
\newcommand{\altsete}[1]{\Lambda^{\exists}_{#1}}
\newcommand{\altsetu}[1]{\Lambda^{\forall}_{#1}}

\newcommand{\vocab}{\gs}

\newcommand{\imply}{\to}
\newcommand{\equi}{\leftrightarrow}
\newcommand{\COOR}{\mbox{\it COOR}}
\newcommand{\VAL}{\mbox{\it VAL}}
\newcommand{\HP}{\mbox{\it HP}}
\newcommand{\ST}{\mbox{\it ST}}
\newcommand{\depth}{\mathop{\sl depth}\nolimits}
\newcommand{\sftree}{q(n;\mbox{\it trees})}
\newcommand{\sftreez}{q_0(n;\mbox{\it trees})}
\newcommand{\sfclass}{q(n;{\cal C})}
\newcommand{\keq}{\mathbin{{\equiv}_k}}
\newcommand{\notkeq}{\mathbin{{\not\equiv}_k}}
\newcommand{\EHRV}[1]{\mbox{\it Ehrv}\,(#1)}

\newcommand{\lpre}[2]{L^{\mbox{\scriptsize\it prenex}}_{#1} (#2)}
\newcommand{\sfpre}[2]{s^{\mbox{\scriptsize\it prenex}}_{#1} (#2)}
\newcommand{\dpre}[1]{D^{\mbox{\scriptsize\it prenex}} (#1)}
\newcommand{\lpree}[1]{L^{\mbox{\scriptsize\it prenex}} (#1)}
\newcommand{\rk}{\mathop{\sl rk}\nolimits}
\newcommand{\sat}[1]{\mbox{\it Sat}\,(#1)}
\newcommand{\finsat}[1]{\mbox{\it Sat\/{}}_{\mbox{\scriptsize\it fin}}(#1)}
\newcommand{\sateq}[1]{\mbox{\it Sat\/{}}^{=}(#1)}
\newcommand{\finsateq}[1]{\mbox{\it Sat\/{}}^{=}_{\mbox{\scriptsize\it fin}}(#1)}
\newcommand{\fintheq}[1]{\mbox{\it Th\/{}}^{=}_{\mbox{\scriptsize\it fin}}(#1)}
\newcommand{\calD}{{\cal D}}
\newcommand{\calS}{{\cal S}}
\newcommand{\calI}{{\cal I}}
\newcommand{\calG}{{\cal G}}
\newcommand{\calC}{{\cal C}}

\title{
Succinct Definitions in\\
the First Order Theory of Graphs}

\author{Oleg Pikhurko\thanks{%
Department of Mathematical Sciences,
Carnegie Mellon University, Pittsburgh, PA 15213-3890.
Web: {\tt http://www.math.cmu.edu/\~{}pikhurko/}}
\quad
Joel Spencer\thanks{%
Courant Institute, New York University, New York, NY 10012.
E-mail: {\tt spencer@cs.nyu.edu}}
\quad
Oleg Verbitsky\thanks{%
Dept.\ of Mechanics \& Mathematics, Kyiv University, Ukraine.
E-mail: {\tt oleg@ov.litech.net}}}

\date{24 March 2004}

\begin{document}
\maketitle

\begin{abstract}
We say that a first order sentence $A$ defines a graph $G$ if $A$ is
true on $G$ but false on any graph non-isomorphic to $G$.
Let $L(G)$ (resp.\ $D(G)$) denote the minimum length (resp.\ quantifier
rank) of a such sentence. We define the succinctness function $s(n)$
(resp.\ its variant $q(n)$) to be the minimum $L(G)$ (resp.\ $D(G)$)
over all graphs on $n$ vertices.

We prove that $s(n)$ and $q(n)$ may be so small that for no general
recursive function $f$ we can have $f(s(n))\ge n$ for all $n$.
However, for the function $q^*(n)=\max_{i\le n}q(i)$, which is
the least monotone nondecreasing function bounding $q(n)$ from above,
we have $q^*(n)=(1+o(1))\log^*n$, where $\log^*n$ equals the minimum
number of iterations of the binary logarithm sufficient to lower $n$
to 1 or below.

We show an upper bound $q(n)<\log^*n+5$ even under the restriction of
the class of graphs to trees. Under this restriction, for $q(n)$
we also have a matching lower bound.

We show a relationship $D(G)\ge(1-o(1))\log^*L(G)$ and prove, using the
upper bound for $q(n)$, that this relationship is tight.

For a non-negative integer $a$, let $D_a(G)$ and $q_a(n)$ denote the
analogs of $D(G)$ and $q(n)$ for defining formulas in the negation normal
form with at most $a$ quantifier alternations in any sequence of nested
quantifiers. We show a superrecursive gap between $D_0(G)$ and $D_3(G)$
and hence between $D_0(G)$ and $D(G)$. Despite it, for $q_0(n)$ we still
have a kind of log-star upper bound: $q_0(n)\le2\log^*n+O(1)$ for
infinitely many~$n$.
\end{abstract}

\clearpage

\mbox{}

\vspace{-30mm}

\mbox{}

\tableofcontents

\clearpage

\section{Introduction}

We study sentences about graphs expressible in the laconic first order
language with two relation symbols $\sim$ and $=$ for, respectively,
the adjacency and the equality relations.
The first order means that we are allowed
to quantify only over vertices, in opposite to second order logic
where we can quantify over sets of vertices. The difference
between the first order and the second order worlds is essential.
In the first order language we cannot express many basic properties
of graphs, as the connectedness, the property of being bipartite etc
(see, e.g., \cite[theorems 2.4.1 and 2.4.2]{Spe}).
On the other hand, the crucial for us fact is that the first order language
is powerful enough to define any individual finite graph up to isomorphism.
Indeed, a graph $G$ with vertex set $V(G)=\{1,\ldots,n\}$
and edge set $E(G)$ is defined by the formula
\begin{equation}\label{eq:defsimpl}
\begin{array}{rccccl}
\exists x_1\ldots\exists x_n\forall x_{n+1}\biggl(&\displaystyle
\bigwedge_{1\le i<j\le n}& \neg (x_i=x_j) &\And
&\displaystyle \bigvee_{i\le n}& x_{n+1}=x_i\\[5mm]
\And&\displaystyle \bigwedge_{\{i,j\}\in E(G)}& x_i\sim x_j &\And&
\displaystyle \bigwedge_{\{i,j\}\notin E(G)}& \neg (x_i\sim x_j)\biggr).
\end{array}
\end{equation}
This fact, though very simple, highlights a fundamental difference
between the finite and the infinite: There are non-isomorphic countable
graphs satisfying precisely the same first order sentences
(see, e.g., \cite[theorem 3.3.2]{Spe}).

The question we address is how succinctly a graph $G$ on $n$ vertices
can be defined by first order means. We consider two natural measures
of succinctness --- the length of a first order formula and its quantifier
rank. The latter is the maximum number of nested quantifiers in the formula.
Let $D(G)$ be the minimum quantifier rank of a closed first order formula
defining $G$, that is, being true on $G$ and false on any
other graph non-isomorphic to $G$. The sentence \refeq{eq:defsimpl}
ensures that $D(G)\le n+1$. This bound generally cannot be improved
as $D(G)=n+1$ for $G$ being the complete or the empty graph on $n$
vertices. However, for all other graphs
we have $D(G)\le n$. Thus, it is reasonable to try to lower the
trivial upper bound of $n+1$ to some $u(n)\le n$ and explicitly
describe all exceptional graphs with $D(G)>u(n)$. This is done in
\cite{PVV} with $u(n)=n/2+O(1)$ (see also \cite{PVe} for a generalization
to arbitrary structures). More precisely, let us
call two vertices of a graph {\em similar\/} if they are simultaneously
adjacent or not to any other vertex. This is an equivalence relation
and each equivalence class spans a complete or an empty subgraph.
Let $\gs(G)$ denote the maximum number of pairwise similar vertices
in $G$. Then, as shown in \cite{PVV},
$$
\gs(G)+1\le D(G)\le\max\ofcurle{\frac{n+5}2,\gs(G)+2}.
$$
It seems doubtful that results of this sort are possible
to obtain with upper bound $u(n)=cn+O(1)$ for each constant $c<1/2$.
The known Cai-F\"urer-Immerman construction \cite{CFI} gives graphs with
linear $D(G)$ that may serve as counterexamples to most natural
conjectures in this direction.

While the paper \cite{PVV} addresses the definability of $n$-vertex graphs
in the worst case, in \cite{KPSV} we treat the average case.
Let $G$ be a random graph distributed uniformly among the graphs
with vertex set $\{1,\ldots,n\}$. Then, as shown in \cite{KPSV},
$$
|D(G) - \log_2n| = O(\log_2\log_2 n)
$$
with probability $1-o(1)$.

We now consider another extremal case of the graph definability
problem. How succinct can be a first order definition of a graph
on $n$ vertices in the best case? Namely, we address the succinctness
function $q(n)$ defined as the minimum
$D(G)$ over $n$-vertex $G$. We also define $L(G)$ to be the minimum
length of a sentence defining $G$ and $s(G)$ to be the minimum $L(G)$
over $n$-vertex $G$. Trivially, $q(n)<s(n)$.
Our first result is that $s(n)$ and $q(n)$ may be so
small that for no general recursive function $f$ we can have
$f(s(n))\ge n$ for all $n$.

The proof is based on simulation of a Turing machine $M$ by a first order
formula $A_M$ in which a computation of $M$
determines a graph satisfying $A_M$ and vise versa.
Such techniques were developed in the classic research on Hilbert's
{\em Entscheidungsproblem\/} by Turing, Trakhtenbrot, B\"uchi and
other researchers (see \cite{BGG} for survey and references).
The novel feature of our simulation is that it works if we restrict
the class of structures to graphs. The key ingredient of our proof
is a gadget allowing us to impose an order relation on the vertex set
of a graph.

As a by-product, we obtain another proof of Lavrov's result \cite{Lav}
that the first order theory of finite graphs is undecidable.
Our proof actually shows the undecidability of the
$\forall^*\exists^p\forall^s\exists^t$-fragment of this theory
for some $p$, $s$, and $t$.

{}From the fact that $q(n)$ and $n$ are not recursively linked,
it easily follows that, if a general recursive function $l(n)$
is monotone nondecreasing and tends to the infinity, then
\begin{equation}\label{eq:intro0}
q(n)<l(n)\mbox{\ \ for infinitely many\ \ }n.
\end{equation}
Our next result establishes a general upper bound
\begin{equation}\label{eq:intro1}
q(n)<\log^*n+5\mbox{\ \ for\ all\ \ }n.
\end{equation}
Here $\log^*n$ equals the minimum number of iterations of the binary
logarithm sufficient to lower $n$ below 1. It turns out that this is
the best possible monotonic upper bound for $q(n)$.
Let $q^*(n)=\max_{i\le n}q(i)$, which is
the least monotone nondecreasing function bounding $q(n)$ from above.
We prove that
\begin{equation}\label{eq:intro2}
q^*(n)\ge\log^*n-\log^*\log^*n-O(1).
\end{equation}
As the upper bound \refeq{eq:intro1} is monotonic, we obtain
\begin{equation}\label{eq:intro3}
q^*(n)=(1+o(1))\log^*n.
\end{equation}
Comparing \refeq{eq:intro3} to \refeq{eq:intro0} with $l(n)=\log^*n$,
we conclude that $q(n)$ infinitely often deviates from its ``smoothed''
version $q^*(n)$ and, in particular, is essentially nonmonotonic.

Proving \refeq{eq:intro1} and \refeq{eq:intro2}, we use a robust
technical tool given by the \EF\/ game \cite{Ehr}
(these techniques were also developed by Fra\"\i ss\'e \cite{Fra}
in a different setting).

As a matter of fact, we prove the upper bound \refeq{eq:intro1} under
the restriction of the class of graphs to trees only, that is, we have
$q(n)\le\sftree<\log^*n+5$. Recall that, by \refeq{eq:intro0}, $q(n)$
is infinitely often so small that we cannot bound it from below by
any ``regular'' function. The proof of this fact cannot be carried through
for $\sftree$ because, as a well-known corollary of the Rabin theorem
\cite{Rab}, the first order theories of both all and finite trees are
decidable and hence a Turing machine computation cannot be simulated
by a first order sentence about trees. In fact, for $\sftree$ we
establish a matching lower bound thereby determining this function
asymptotically, namely,
$$
\sftree=(1+o(1))\log^*n.
$$

We pay a special attention to defining sentences having a restricted
structure.
For a non-negative integer $a$, let $D_a(G)$ and $q_a(n)$ denote the
analogs of $D(G)$ and $q(n)$ for defining formulas in the negation normal
form with at most $a$ quantifier alternations in any sequence of nested
quantifiers. The superrecursive gap between $s(n)$ and $n$ is actually
shown even under the restriction of the alternation number to 3.
Note also that, as follows from a result in \cite{KPSV},
$q_3(n)\le\log^*n+O(1)$ and hence \refeq{eq:intro3} holds with
alternation number~3.

On the other hand,
we show a superrecursive gap between $D_0(G)$ and $D_3(G)$
and hence between $D_0(G)$ and $D(G)$. Despite it, for $q_0(n)$ we also
have a kind of log-star upper bound: $q_0(n)\le2\log^*n+O(1)$ for
infinitely many~$n$. It is worth noting that this is not the first case
that we have close results for the alternation number 0 and for the unbounded
alternation number. In \cite{KPSV} we prove that for a random graph
$D(G)$ and $D_0(G)$ are not so far apart from each other. Namely,
$D_0(G)\le(2+o(1))\log_2n$ with probability $1-o(1)$. Yet another
result showing the same phenomenon is obtained in \cite{PVV}.
Given non-isomorphic graphs $G$ and $G'$, let $D(G,G')$ (resp.\ $D_0(G,G')$)
denote the minimum quantifier rank of a sentence (resp.\ in the negation
normal form with no quantifier alternation) which is true on exactly one
of the graphs. As shown in \cite{PVV}, if both $G$ and $G'$ have
$n$ vertices, then $D(G,G')\le D_0(G,G')\le (n+5)/2$ and
there are simple examples of such $G$ and $G'$ with $D(G,G')\ge(n+1)/2$.
Note that logical distinguishing non-isomorphic graphs with equal number
of vertices has close connections to graph canonization algorithms
(see, e.g., \cite{CFI,Gro,PVV} and a monograph~\cite{Imm}).

Relating $D(G)$ and $L(G)$ to one another, we show that
$$
D(G)\ge(1-o(1))\log^*L(G).
$$
Using the bound \refeq{eq:intro1}, we show that this relationship is tight.

Focusing on defining formulas of restricted structure, we also
consider prenex formulas. A superrecursive gap between
$s(n)$ and $n$ can actually be shown under the restriction to this
class. Nevertheless, prenex formulas generally are not competetive
against defining formulas with no restriction on structure.
We observe that graphs showing a huge gap between $D(G)$ and $L(G)$
at the same time show a huge gap between $D(G)$ and its version for
prenex defining formulas.

In conclusion, note that all of our results carry over to general structures
over any relational vocabulary with at least one non-unary relation
symbol. For the upper bounds this claim is straightforward because
graphs can be viewed as a subclass of such structures which is
distinguishable by a single first order sentence. The lower bounds hold
true with minor changes in the proofs.

\section{Background}\label{s:back}

\subsection{Arithmetics}

We define the {\em tower function\/} $T(i)$ by $T(0)=1$ and
$T(i)=2^{T(i-1)}$ for each subsequent $i$.
Sometimes this function will be denoted by $\mbox{\it Tower\,}(i)$.
Given a function $f$, by $f^{(i)}$ we will denote the
$i$-fold composition of $f$. In particular, $f^{(0)}(x)=x$.
By $\log n$ we always mean the logarithm base 2. The inverse of the tower
function, the {\em log-star\/} function $\log^*n$, is defined
by $\log^*n=\min\setdef{i}{T(i)\ge n}$. For a real $x$, the notation
$\lceil x\rceil$ (resp.\ $\lfloor x\rfloor$) stands for the integer nearest
to $x$ from above (resp.\ from below).

\subsection{Graphs}\label{ss:graphs}

Given a graph $G$, we denote its vertex set by $V(G)$ and its edge set
by $E(G)$. The {\em order\/} of $G$, the number of vertices of $G$,
will be sometimes denoted by $|G|$, that is,
$|G|=|V(G)|$. The {\em neighborhood\/} of a vertex $v$ consists of
all vertices adjacent to $v$.
A set $S\subseteq V(G)$ is called {\em independent\/}
if it contains no pair of adjacent vertices.
If $X\subseteq V(G)$, then $G[X]$ denotes the subgraph {\em induced\/}
by $G$ on $X$ (or {\em spanned\/} by $X$ in $G$).
If $u\in V(G)$, then $G-u=G[V(G)\setminus\{u\}]$ is the result
of removal from $G$ the vertex $u$ along with all incident edges.

The {\em distance\/} between vertices $u$ and $v$, the minimum
length of a path connecting the two vertices, is denoted by $d(u,v)$.
If $u$ and $v$ are in different connected components of a graph, then
$d(u,v)=\infty$. The {\em eccentricity\/} of a vertex $v$ is defined
by $e(v)=\max_{u\in V(G)}d(v,u)$. The {\em diameter\/} and the
{\em radius\/} of a graph $G$ are defined by $d(G)=\max_{v\in V(G)}e(v)$
and $r(G)=\min_{v\in V(G)}e(v)$ respectively.
A path in a graph is {\em diametral\/} if its length is equal to
the diameter of the graph. A vertex $v$ is {\em central\/} if $e(v)=r(G)$.

\begin{proposition}\label{prop:ore}
{\bf \cite[Theorem 4.2.2]{Ore}}
Let $T$ be a tree. If $d(T)$ is even, then $T$ has a unique central
vertex $c$ and all diametral paths go through $c$.
If $d(T)$ is odd, then $T$ has exactly two central vertices $c_1$ and
$c_2$ and all diametral paths go through the edge $\{c_1,c_2\}$.
\end{proposition}

\subsection{Logic}\label{ss:logic}

\subsubsection{Formulas}

First order formulas are assumed to be over the set of connectives
$\{\neg,\And,\Or\}$.
A {\em sequence of quantifiers\/} is a finite word over the alphabet
$\{\exists,\forall\}$. If $S$ is a set of such sequences, then
$\exists S$ (resp.\ $\forall S$) means the set of concatenations
$\exists s$ (resp.\ $\forall s$) for all $s\in S$. If $s$ is a sequence
of quantifiers, then $\bar s$ denotes the result of replacement of all
occurrences of $\exists$ to $\forall$ and vise versa in $s$. The
set $\bar S$ consists of all $\bar s$ for $s\in S$.

Given a first order formula ${A}$, its set of {\em sequences of nested
quantifiers\/} is denoted by $\nest{A}$ and defined by induction as
follows:
\begin{enumerate}
\item
$\nest{A}=\{ \epsilon \}$ if ${A}$ is atomic; here $\epsilon$ denotes the
empty word.
\item
$\nest{\neg{A}}=\overline{\nest{A}}$.
\item
$\nest{{A}\And{B}}=\nest{{A}\Or{B}}=\nest{A}\cup\nest{B}$.
\item
$\nest{\exists x{A}}=\exists\nest{A}$ and
$\nest{\forall x{A}}=\forall\nest{A}$.
\end{enumerate}
The {\em quantifier rank\/} of a formula ${A}$, denoted by $\qr{A}$
is the maximum length of a string in $\nest{A}$.

We adopt the notion of the {\em alternation number\/} of a formula
(cf.\ \cite[Definition 2.8]{Pez}).
Given a sequence of quantifiers $s$, let $\alt s$ denote the number
of occurrences of $\exists\forall$ and $\forall\exists$ in $s$.
The {\em alternation number\/} of a first order formula ${A}$,
denoted by $\alt{A}$, is the maximum $\alt s$ over $s\in\nest{A}$.
The alternation number has an absolutely clear meaning for
formulas in the {\em negation normal form}, where
the connective $\neg$ occurs only in front of atomic subformulas.
This number is defined for any formula $A$ so that, if $A$ is reduced
to an equivalent formula $A'$ in the negation normal form, then
$\alt A=\alt{A'}$.

Viewing a formula $A$ as a string of symbols over
the countable first order alphabet (where each variable
and each relation is denoted by a single symbol), we denote the
{\em length\/} of $A$ by $|A|$. Note that if one prefers, in a natural way,
to encode variable and relation symbols in a finite alphabet,
then the length will increase but stay within $|A|\log|A|$.

We call $A$ an {\em $\exists$-formula\/} (resp.\ {\em $\forall$-formula\/})
if any sequence in $\nest A$ with maximum number of quantifier
alternations starts with $\exists$ (resp.\ $\forall$).
We denote the set of formulas in the negation normal form
with alternation number at most $m$ by $\altset m$.
By $\altsete m$ (resp.\ $\altsetu m$) we denote
the subset of $\altset m$ consisting of formulas in $\altset{m-1}$
and $\exists$-formulas (resp.\ $\forall$-formulas) in
$\altset m\setminus\altset{m-1}$.
We will call formulas in $\altsete 0$ and $\altsetu 0$
{\em existential\/} and {\em universal\/} respectively.

A {\em prenex formula\/} is a formula with all its quantifiers up front.
In this case there is a single sequence of nested quantifiers and
the quantifier rank is just the number of quantifiers
occurring in a formula.
Let $\Sigma_1$ and $\Pi_1$ denote, respectively, the sets of
existential and universal prenex formulas. Furthermore,
let $\Sigma_m$ (resp.\ $\Pi_m$) be the extension of
$\Sigma_{m-1}\cup\Pi_{m-1}$ with prenex formulas in
$\altsete{m-1}$ (resp.\ $\altsetu{m-1}$).
Note that the classes of formulas $\altset m$, $\altsete{m}$, $\altsetu{m}$,
$\Sigma_m$, and $\Pi_m$ are defined so that they are closed with respect to
subformulas.

The following lemma is an immediate consequence of the standard
reduction of a formula to the prenex form.

\begin{lemma}\label{lem:conjtoprenex}
The conjunction of $\Sigma_m$-formulas (resp.\ $\Pi_m$-formulas) is
effectively reducible to an equivalent $\Sigma_m$-formula
(resp.\ $\Pi_m$-formula). The same holds for the disjunction.\noproof
\end{lemma}

We write $A\equiv B$ if $A$ and $B$ are logically equivalent formulas
and $A\doteq B$ if $A$ and $B$ are literally the same.

\begin{lemma}\label{lem:alttoprenex}
\mbox{}

\begin{enumerate}
\item
Any formula in $\altsete m$ is effectively reducible to an equivalent
formula in $\Sigma_{m+1}$.
\item
Any formula in $\altsetu m$ is effectively reducible to an equivalent
formula in $\Pi_{m+1}$.
\item
Any formula in $\altset m$ is effectively reducible to an equivalent
formula in $\Sigma_{m+2}$ or, as well, to an equivalent
formula in $\Pi_{m+2}$.
\end{enumerate}
\end{lemma}

\begin{proof}
Item 3 follows from Items 1 and 2 as $\altset m$ is included both in
$\altsete{m+1}$ and $\altsetu{m+1}$. To prove Items 1 and 2, we proceed
by induction on $m$.

Consider the base case of $m=0$.
Assume that $A\in\altsete 0$ and let $t=t(A)$ denote the total number of
quantifiers and connectives $\And$, $\Or$ in $A$. We prove that
$A$ has an equivalent formula $A'\in\Sigma_1$ using induction on~$t$.
If $t=0$, then $A$ is quantifier-free and hence in $\Sigma_0$.
Let $t\ge 1$. Assume that $A\doteq\exists x B$.
Since $t(B)=t(A)-1$, the assumption of induction on $t$ applies to $B$.
Therefore $B$ reduces to an equivalent formula $B'\in\Sigma_1$
and we set $A'=\exists x B'$. Assume that $A\doteq B\And C$
(the case that $A\doteq B\Or C$ is similar). Neither of $t(B)$
and $t(C)$ exceeds $t(A)-1$ and, by the assumption of induction on $t$,
for $B$ and $C$ we have equivalents $B'$ and $C'$ in $\Sigma_1$. Then
$A\equiv B'\And C'$ reduces to an equivalent in $\Sigma_1$ by
Lemma~\ref{lem:conjtoprenex}.

The reducibility of $\altsetu 0$ to $\Pi_1$ is proved similarly.

Let $m\ge 1$ and assume that Items 1 and 2 of the lemma are true
for the preceding value of $m$. Given $A\in\altsete m$, we show how to find
an equivalent formula $A'\in\Sigma_{m+1}$
(the reduction of $\altsetu m$ to $\Pi_{m+1}$ is similar).
We again use induction on $t=t(A)$.
If $t=0$, then $A$ is in $\Sigma_0$.
Let $t\ge 1$. If $A\doteq\forall x B$, then $A\in\altsetu{m-1}$
and, by the assumption of induction on $m$, $A$ has an equivalent
$A'\in\Pi_m\subset\Sigma_{m+1}$.
If $A\doteq\exists x B$, $A\doteq B\And C$, or $A\doteq B\Or C$,
then $B,C\in\altsete m$ and both $t(B)$ and $t(C)$ are smaller than $t(A)$.
We are done by the assumption of induction on $t$ and
Lemma~\ref{lem:conjtoprenex}.
\end{proof}

A formula with all variables bound is called a {\em closed formula\/} or
a {\em sentence}.

\begin{lemma}\label{lem:lvsdprenex}
If $A$ is a closed prenex formula of quantifier rank $q$ with occurrences
of $h$ binary relation symbols, then it can be rewritten in an equivalent
form $A'$ with the same quantifier prefix so that $|A'|=O(hq^22^{hq^2})$.
\end{lemma}

\begin{proof}
Let $B(x_1,\ldots,x_q)$ be the quantifier-free part of $A$.
The $B$ is a Boolean combination of $m=h{q\choose 2}$ atomic subformulas
and hence is representible as a DNF of length $O(m2^m)$.
\end{proof}

\subsubsection{Structures}

A {\em relational vocabulary\/} $\vocab$ is a finite set of {\em relation
symbols\/} augmented with their {\em arities}.
We always assume the presence of the binary relation symbol $=$
standing for the equality relation and do not include it in~$\vocab$.
The only exception will be Subsection \ref{ss:tfg} where the presence
or the absence of equality will be stated explicitly.

A {\em structure over vocabulary
$\vocab$\/} (or an $\vocab$-structure) is a set along with relations
that are named by symbols in $\vocab$ and have the corresponding
arities. We mostly deal with the vocabulary of a single
binary relation symbol. A structure over this vocabulary can be viewed
as a {\em directed graph\/} (or {\em digraph\/}).
We treat {\em graphs\/} as structures with a single binary relation which
is symmetric and anti-relexive. This relation will be called the
{\em adjacency\/} relation and denoted by $\sim$.

If all relation symbols of a sentence $A$ are
from the vocabulary $\vocab$ and $G$ is an $\vocab$-structure, then
$A$ is either true or false on $G$. In the former case $G$ is called
a {\em model\/} of $A$. We also say that $G$ {\em satisfies\/} $A$.
We call $A$ {\em valid\/} if all $\vocab$-structures satisfy $A$.
We call $A$ {\em (finitely) satisfiable\/} if it has a (finite) model.
Clearly, $A$ is valid iff $\neg A$ is unsatisfiable.

\subsubsection{Computability}

Whenever we say that something can be done {\em effectively}, we
mean that this can be implemented by an {\em algorithm}.
No restrictions on running time or space are assumed.
Professing {\em Church's thesis}, we here do not specify any
definition of the algorithm. Nevertheless, we will refer to
{\em Turing machines\/} (see Subsection \ref{ss:tmdef}) and
{\em recursive functions\/} in Sections \ref{s:tm} and \ref{s:other}.
As a basic fact, these two computational models are
equally powerful, under an effective bijection between binary words
and non-negative integer numbers.

Let $X$ be a set of words over a finite alphabet.
The {\em decision problem\/} for $X$ is the problem of recognition
whether or not a given word belongs to $X$. If there is an algorithm
that does it, the decision problem is {\em solvable}
(or $X$ is {\em decidable}).

The {\em halting problem\/} is the problem of deciding, for given
Turing machine $M$ and input word $w$, whether $M$ eventually halts
on $w$ or runs forever. This is a basic unsolvable problem.
It is well known that, if we fix $w$ to be the empty word, the restricted
problem remains unsolvable.

The {\em (finite) satisfiability problem\/}
is the problem of recognizing whether or not a
given sentence is (finitely) satisfiable
(we here assume any natural encoding of formulas in a finite alphabet).
Settling Hilbert's {\em Entscheidungsproblem}, Church and Turing proved
that the satisfiability problem is unsolvable. The unsolvability of
the finite satisfiability problem was shown by Trakhtenbrot~\cite{Tra}.

A {\em general recursive function\/} is an everywhere defined recursive
function.

\subsubsection{The \BS\/ class of formulas and the Ramsey theorem}

A class of formulas has the {\em finite model property\/} if
every satisfiable formula in the class has a finite model.
By the completeness of the predicate calculus with equality,
the set of valid sentences is recursively enumerable.
{}From here it is not hard to conclude that,
if a class of formulas has the finite model property,
the satisfiability and the finite satisfiability problems
for this class are solvable.

The {\em\BS\/ class\/} consists of prenex formulas in which the existential
quantifiers all precede the universal quantifiers, that is, this is
another name for~$\Sigma_2$.

\begin{proposition}
{\bf (The Ramsey theorem \cite{Ram}\footnote{%
The {\em combinatorial Ramsey theorem}, a cornerstone of
{\em Ramsey theory}, appeared in this paper as a technical tool.%
})}
For each vocabulary $\vocab$ there is a
general recursive function $f\function{\bf N}{\bf N}$
such that the following is true:
Assume that a $\vocab$-sentence $A$ with equality is in the \BS\/ class.
If $A$ has a model of some cardinality at least $f(\qr A)$
(possibly infinite), then it has a model in every cardinality at least
$f(\qr A)$. As a consequence, the \BS\/ class
of formulas with equality has the finite model property and
hence both the satisfiability and the finite satisfiability problems
restricted to this class are solvable.
\end{proposition}

\subsubsection{Definability}

Let $G$ and $G'$ be non-isomorphic graphs
and ${A}$ be a first order sentence with equality over vocabulary
$\{{\sim}\}$.
We say that ${A}$ {\em distinguishes $G$ from $G'$} if
${A}$ is true on $G$ but false on $G'$.
By $D(G,G')$ (resp.\ $D_k(G,G')$) we denote the minimum quantifier rank of
a sentence (with alternation number at most $k$ resp.) distinguishing $G$
from $G'$.

We say that a sentence ${A}$ {\em defines\/} a graph $G$ (up to isomorphism)
if ${A}$ distinguishes $G$ from any non-isomorphic graph $G'$.
To ensure that $A$ has no other models except graphs,
we will tacitly assume that $A$ has form $A\doteq
\forall_x(x\not\sim x\And\forall_y(x\sim y\imply y\sim x))\And B$.
By $D(G)$ (resp.\ $D_a(G)$) we denote the minimum
quantifier rank of a sentence defining $G$
(with alternation number at most $a$ resp.).
By $L(G)$ (resp.\ $L_a(G)$) we denote the minimum
length of a sentence defining $G$
(with alternation number at most $a$ resp.).

A sentence is called {\em defining\/} if it defines a graph.
Note that any defining sentence must contain the equality symbol.
Let us stress that graphs $G'$ in the above definition may have any
cardinality.

\begin{lemma}\label{lem:findef}
All finite graphs and only finite graphs posses defining sentences.
\end{lemma}

\begin{proof}
Any finite graph is indeed definable as it has at least
the wasteful definition~\refeq{eq:defsimpl}.
By the Upward L\"owenheim-Skolem theorem (see \cite[corollary 2.35]{Men}),
if a sentence with equality has an infinite model, it has a model
of any infinite cardinality. By this reason, no infinite graph
has defining sentence in the sence of our definition.
\end{proof}

\begin{lemma}\label{lem:defdec}
The class of defining $\altsete 1$-sentences is decidable.
\end{lemma}

\begin{proof}
Suppose that we are given a sentence $A\in\altsete1$.
By Lemma \ref{lem:alttoprenex} (1), we can reduce it to an equivalent
formula in the \BS\/ class and apply the Ramsey theorem.
We are able to recognize if $A$ is defining in four steps.
\begin{enumerate}
\item
Check if $A$ is finitely satisfiable.
\item
If so, trying graphs one by one, we eventually find a graph of the
smallest order $n$ satisfying $A$ (this is actually done in the first step,
if it is based directly on the Ramsey theorem).
\item
Check if there is any other graph of order $n$ satisfying~$A$.
\item
If not, check if a $\altsete1$-sentence
$A\And\exists_{x_1,\ldots,x_{n+1}}(\bigwedge_{1\le i<j\le n+1}x_i\ne x_j)$
is satisfiable.
\end{enumerate}
If not, and only in this case, $A$ is defining.
\end{proof}

\section{The \EF\/ game}\label{s:game}

In this section we borrow a lot of material from \cite[section 2]{Spe}.
To make our exposition self-contained, we sketch some proofs that
can be found in \cite{Spe} in more detail.

The {\em \EF\/ game\/} is played on a pair of structures of the same
vocabulary. We give the definition conformable to the case of graphs.

Let $G$ and $H$ be graphs with disjoint vertex sets.
The $k$-round \EF\/ game on $G$ and $H$,
denoted by $\game_k(G,H)$, is played by
two players, Spoiler and Duplicator, with $k$ pairwise distinct
pebbles $p_1,\ldots,p_k$, each given in duplicate. Spoiler starts the game.
A {\em round\/} consists of a move of Spoiler followed by a move of
Duplicator. At the $i$-th move Spoiler takes pebble $p_i$, selects one
of the graphs $G$ or $H$, and places $p_i$ on a vertex of this graph.
In response Duplicator should place the other copy of $p_i$ on a vertex
of the other graph. It is allowed to
place more than one pebble on the same vertex.

Let $u_i$ (resp.\ $v_i$)
denote the vertex of $G$ (resp.\ $H$) occupied by $p_i$, irrespectively
of who of the players placed the pebble on this vertex.
If
$$
u_i=u_j\mbox{ iff } v_i=v_j\mbox{ for all }1\le i<j\le k,
$$
and the component-wise correspondence $(u_1,\ldots,u_k)$ to
$(v_1,\ldots,v_k)$ is a partial isomorphism from $G$ to $H$, this is
a win for Duplicator;  Otherwise the winner is Spoiler.

The {\em $a$-alternation\/} \EF\/ game on $G$ and $H$ is a variant
of the game in which Spoiler is allowed to switch from one graph
to another at most $a$ times during the game, i.e., in at most $a$
rounds he can choose the graph other than that in the preceding round.

Let $0\le s\le k$, $r=k-s$, and assume that at the start of the game
the pebbles $p_1,\ldots,p_s$ are already on the board
at vertices $\bar u=u_1,\ldots,u_s$ of $G$ and $\bar v=v_1,\ldots,v_s$
of $H$. The $r$-round game with this initial configuration
is denoted by $\game_r(G,\bar u,H,\bar v)$.
We write $G,\bar u\keq H,\bar v$ if Duplicator has a winning strategy
in this game.

It is not hard to check that $\keq$ is an equivalence relation.
The $k$-Ehrenfeucht value
of a graph $G$ with vertices $u_1,\ldots,u_s$ marked by pebbles is
the equivalence class it belongs to under $\keq$.  We let
$\EHRV{k,s}$ denote the set of all possible $k$-Ehrenfeucht values
for graphs with $s$ marked vertices. Let $\EHRV k=\EHRV{k,0}$
denote the set of $k$-Ehrenfeucht values for graphs (with no marked vertex).

\begin{lemma}\label{lem:nextmove}
Assume that $s<k$. Let $\bar u=u_1,\ldots,u_s$ and
$S(G,\bar u)$ denote the
set of $\keq$-equivalence classes of $G$ with $s+1$ marked vertices
$\bar u,u$ for all $u\in G\setminus\{u_1,\ldots,u_s\}$.
Then $G,\bar u\keq H,\bar v$ iff $S(G,\bar u)=S(H,\bar v)$.
\end{lemma}

\begin{proof}
Consider the game $\game_{k-s}(G,\bar u,H,\bar v)$.
Suppose that $S(G,\bar u)\ne S(H,\bar v)$, for example, there is
$u\in V(G)$ such that $G,\bar u,u\notkeq H,\bar v,v$
for any $v\in V(H)$.
Let Spoiler select this $u$ and let $v$ denote Duplicator's response.
{}From now on the players actually play $\game_{k-s-1}(G,\bar u,u,H,\bar v,v)$,
where Spoiler has a winning strategy.

Suppose that $S(G,\bar u)=S(H,\bar v)$.
If Spoiler selects, for example, a vertex $u\in V(G)$, then
Duplicator responds with $v\in V(H)$ such that $G,\bar u,u\keq H,\bar v,v$
and hence has a winning strategy in the remaining part of the game.
\end{proof}

\begin{lemma}\label{lem:ehrvks}
{\bf \cite[theorem 2.2.1]{Spe}}
For any $s$ and $k$, $\EHRV{k,s}$ is a finite set. Furthermore,
let $f(k,s)=|\EHRV{k,s}|$. Then
\begin{eqnarray}
f(k,k)&\le&4^{k\choose2},\label{eq:kk}\\
f(k,s)&\le&2^{f(k,s+1)}\label{eq:ks}
\end{eqnarray}
for $s<k$.
\end{lemma}

\begin{proof}
The bound \refeq{eq:kk} holds because the $\keq$-equivalence class
of $G$ with marked $u_1,\ldots,u_k$ is determined by the equality relation
on the sequence $u_1,\ldots,u_k$ and the induced subgraph
$G[\{u_1,\ldots,u_k\}]$.
The bound \refeq{eq:ks} holds because the $\keq$-equivalence class
of an arbitrary $G$ with marked $\bar u=u_1,\ldots,u_s$ is, according to
Lemma \ref{lem:nextmove}, determined by $S(G,\bar u)$, a subset of
$\EHRV{k,s+1}$.
\end{proof}

As a consequence, we obtain the following bound.

\begin{lemma}\label{lem:ehrvk}
{\bf \cite[theorem 2.2.2]{Spe}}
$|\EHRV k|\le T(k+2+\log^*k)+O(1)$.
\noproof
\end{lemma}

We say that a formula $A(x_1,\ldots,x_s)$ with $s$ free variables
{\em defines\/} an Ehrenfeucht value $\ga\in\EHRV{k,s}$ if
$A$ is true on a graph $G$ with variables $x_1,\ldots,x_s$ assigned
vertices $u_1,\ldots,u_s$ for exactly those $G,u_1,\ldots,u_s$
which are in~$\ga$.

\begin{lemma}\label{lem:defehrv}
{\bf \cite[theorem 2.3.2]{Spe}}
For any $\ga\in\EHRV{k,s}$ there is a formula $A_\ga$
with $\qr{A_\ga}=k-s$ that defines $\ga$. Moreover,
\begin{eqnarray}
|A_\ga|&\le&18{k\choose 2}\mbox{\ if\ }s=k\mbox{\ and\ }
\label{eq:akk}\\
|A_\ga|&\le&f(k,s+1)\of{\max\setdef{|A_\gb|}{\gb\in\EHRV{k,s+1}}+10}
\mbox{\ if\ }s<k.\label{eq:aks}
\end{eqnarray}
\end{lemma}

\begin{proof}
The bound \refeq{eq:akk} holds because every $\ga\in\EHRV{k,k}$
is defined by a formula of the type
$$
\bigwedge_{1\le i<j\le k}(*(x_i=x_j)\wedge\star(x_i\sim x_j)),
$$
where $*$ and $\star$ is $\neg$ for some of $(i,j)$ and nothing for the others,
depending on adjacencies among the marked vertices of a
$G,u_1,\ldots,u_k$ in~$\ga$.

Let $s<k$ and assume that every $\beta\in \EHRV{k,s+1}$
has a defining formula $A_{\beta}(x_1,\ldots,x_s,x)$ of quantifier
rank $k-s-1$. Consider an $\ga\in \EHRV{k,s}$ and choose a representative
$G,\bar u$ of $\ga$. Define $S(\ga)=S(G,\bar u)$, where the right hand
side is as in Lemma \ref{lem:nextmove}. By this lemma, the definition
does not depend on a particular choice of $G,\bar u$. We set
$$
A_{\ga}(x_1,\ldots,x_s)\doteq
\bigwedge_{\beta\in S(\ga)}\exists_x
A_{\beta}(x_1,\ldots,x_s,x) \wedge \bigwedge_{\beta\notin S(\ga)}\neg\exists_x
A_{\beta}(x_1,\ldots,x_s,x).
$$
It is clear that $G$ with designated $\bar u=u_1,\ldots,u_s$ satisfies
$A_{\ga}$ iff the set of Ehrenfeucht values with additional designated $u$
is equal to $S(\ah)$. By Lemma \ref{lem:nextmove}, the latter condition is
true iff $G,\bar u$ has Ehrenfeucht value~$\ga$.
\end{proof}

\begin{proposition}\label{prop:loggames}
Suppose that $G$ and $H$ are non-isomorphic graphs.

\begin{enumerate}
\item
Let $R(G,H)$ denote the minimum $k$ such that $G$ and $H$ have
different $k$-Ehrenfeucht values. Then $D(G,H)=R(G,H)$. In other words,
$D(G,H)$ equals the minimum $k$ such that Spoiler has a winning
strategy in $\game_k(G,H)$.
\item
$D_a(G,H)$ equals the minimum $k$ such that Spoiler has a winning
strategy in the $a$-alternation $\game_k(G,H)$.
\end{enumerate}
\end{proposition}

\noindent
We refer the reader to \cite[Theorem 2.3.1]{Spe}
for the proof of the first claim and to \cite{Pez} for the second claim.

\begin{proposition}\label{prop:ddd}
\begin{eqnarray*}
D(G)&=&\max\setdef{D(G,H)}{H\mbox{\ and\ } G\mbox{\ are\ non-isomorphic}},\\
D_a(G)&=&\max\setdef{D_a(G,H)}{H\mbox{\ and\ } G\mbox{\ are\ non-isomorphic}}.
\end{eqnarray*}
The first equality can be restated as follows: $D(G)$ equals the minimum
$k$ such that the $k$-Ehrenfeucht value of $G$ contains only graphs
isomorphic to~$G$.
\end{proposition}

\begin{proof}
We give a proof of the first equality that can be easily adopted for
the second equality.
Denote the maximum in the right hand side by $k$.
We have $k\le D(G)$ as a matter of definition.
Conversely, let $\ga\in\EHRV k$ be the class containing $G$.
By Proposition \ref{prop:loggames}, $G$ is, up to isomorphism,
the only member of $\ga$.
For each $\gb\ne\ga$ in $\EHRV k$, fix a representative $H_\gb$.
Let $C_\gb$ be a sentence of quantifier rank at most $k$ distinguishing
$G$ from $H_\gb$. We use Lemma \ref{lem:ehrvks} saying that
$\EHRV k$ is finite. The conjunction of all $C_\gb$ defines $G$
and has quantifier rank $k$. Thus, $D(G)\le k$.
(Alternatively, we could use the known fact that, over a finite vocabulary,
there are only finitely many inequivalent sentences
of bounded quantifier rank, cf.\ Lemma \ref{lem:uniset}.)
\end{proof}

\section{A superrecursive gap: Simulating a Turing machine}\label{s:tm}

\begin{definition}
We define the {\em succinctness function\/} $s(n)$ (for formula length) by
$$
s(n)=\min_{|G|=n} L(G).
$$
The variants with bounded alternation number are defined by
$$
s_a(n)=\min_{|G|=n} L_a(G)
$$
for each $a\ge0$.
\end{definition}

It turns out that $s(n)$ can be so small with respect to $n$ that
the gap between the two numbers cannot be bounded by any recursive
function.

\begin{theorem}\label{thm:hugegap}
There is no general recursive function $f$ such that
\begin{equation}\label{eq:fsn}
f(s_3(n))\ge n\mbox{\ \ for\ all\ \ }n.
\end{equation}
\end{theorem}

\begin{lemma}\label{lem:tmsim}
{\bf (Simulation Lemma)}
Given a Turing machine $M$ with $k$ states, one can effectively
construct a sentence $A_M$ with single binary relation symbol $\sim$
and equality so that the following conditions are met.
\begin{enumerate}
\item
$\qr {A_M}=k+16$.
\item
$|A_M|=O(k^2)$.
\item
$\alt {A_M}=3$.
\item
$A_M$ is effectively reducible to an equivalent prenex formula $P_M$ whose
quantifier prefix has length $k+O(1)$, begins with $k$ existential
quantifiers, and has 3 quantifier alternations.
\item
Any model of $A_M$ is a graph.
If $M$ halts on the empty input word,
then $A_M$ has a unique model $G_M$ and the order of $G_M$ is bigger
than the running time of $M$.
\item
$M$ halts on the empty input word iff $A_M$ has a finite model.
\end{enumerate}
\end{lemma}

\begin{proofof}{Theorem~\ref{thm:hugegap}}
Let $g(k)$ denote the longest running time on the empty input word
$\epsilon$ of a $k$-state Turing machine (non-halting machines
are excluded from consideration).
Recognition whether or not a given Turing machine with $k$ states
halts on $\epsilon$ easily reduces to computation of $g(k)$.
As this variant of the halting problem is well known to be undecidable,
the function $g(k)$ cannot be bounded from above by any general recursive
function. For each $k$, fix a machine $M_k$ with $k$ states whose
running time attains
$g(k)$. Let $A_{M_k}$ be as in the Simulation Lemma, $G_k$ be the model
of $A_{M_k}$, and $n_k$ be the order of $G_k$. Let $l(k)=ck^2$
be the upper bound for $|A_{M_k}|$ ensured by the lemma.
Note that $A_{M_k}$ defines~$G_k$.

Suppose on the contrary that \refeq{eq:fsn} is true for some general
recursive $f$. Since $s_3(n_k)\le|A_{M_k}|\le l(k)$, for every $k$
we have
$$
g(k)<n_k\le f(s_3(n_k))\le\max_{i\le l(k)} f(i),
$$
a contradiction.
\end{proofof}

The proof of the Simulation Lemma takes the rest of this section.

\subsection{Gadgets}

We enrich our language with connectives $\imply$ and $\equi$ for the
implication and the equivalence. Since the alternation number was
defined for formulas with connectives $\neg,\And,\Or$, we should
stress that $\imply$ and $\equi$ are used as shorthands for their
standard definitions through $\neg,\And,\Or$.
We introduce the new uniqueness quantifier
$\exists!$ by
$$
\exists!_xF(x)\doteq\exists_xF(x)\And
\forall_x\forall_y(F(x)\And F(y)\imply x=y)
$$
for any formula $F$ with a free variable $x$ and with no free
occurrences of $y$. Note that one occurrence of the
uniqueness quantifier contributes 2 in the quantifier rank
and 1 in the alternation number.
We use relativized versions of the existential and the universal
quantifiers in the standard way:
\begin{eqnarray*}
\exists_{C(x)}F(x)&\doteq&\exists_x(C(x)\And F(x)),\\
\forall_{C(x)}F(x)&\doteq&\forall_x(C(x)\imply F(x)).
\end{eqnarray*}

To ensure that any model of $A_M$ is a graph, we put in $A_M$
the two graph axioms (the irreflexivity and the symmetry of
the relation $\sim$).

\subsubsection{Ordering}\label{ss:order}

We give a formula $P(x,x')$ with two free variables $x$ and $x'$
that, in any model, shall determine an order on the neighborhood
of $x$. Let $X=\{y:y\sim x\}$ and $X'=\{z:z\sim x'\}$.
Then $P(x,x')$ is the conjunction of the following:
\begin{description}
\item[(P1)]
$\{x,x'\},X,X'$ are all disjoint and each of them is independent.
\item[(P2)]
$\forall_{y\in X}\exists_{z\in X'}\,\,y\sim z$
\item[(P3)]
$\exists_{y\in X}\exists!_{z\in X'}\,\,y\sim z$
\item[(P4)]
$\exists_{y\in X}\forall_{z\in X'}\,\,y\sim z$
\item[(P5)]
$\forall_{y_1\in X}\forall_{y_2\in X}[
\forall_{z\in X'} (y_1\sim z\ra y_2\sim z)\vee
\forall_{z\in X'} (y_2\sim z\ra y_1\sim z)]$
\item[(P6)]
$\forall_{y_1\in X}\forall_{y_2\in X}
[y_1\neq y_2 \ra
\exists_{z\in X'}(y_1\sim z \equi y_2\not\sim z)]$
\item[(P7)]
$\forall_{y\in X}[\exists_{z\in X'}\,\,y\not\sim z \ra
\exists_{y^+\in X}\exists!_{z\in X'} (y^+\sim z \wedge y\not\sim z)]$
\item[(P8)]
$\forall_{y\in X}[\exists!_{z\in X'}\,\,y\sim z\vee
\exists_{y^-\in X}\exists!_{z\in X'} (y\sim z \wedge y^-\not\sim z)]$
\end{description}
Note that $\qr P=4$, $\alt P=2$ (contributed by (P7) and (P8)),
and $|P|=O(1)$.

Consider finite models of $P(x,x')$.  For $y\in X$ let $N^*(y)$ be those
$z\in X'$ adjacent to $y$.  The $N^*(y)$ are distinct (P6), linearly ordered
under inclusion (P5), are nonempty (P2), include a singleton (P3) and all of
$X'$ (P4), and the set of all cardinalities $|N^*(y)|$ has no gaps
(either (P7) or (P8)).  So we must have $|X|=|X'|$ and the elements
can be ordered $x_1,\ldots,x_s$, $x_1',\ldots,x_s'$ so that $x_i,x_j'$ are
adjacent precisely when $j\leq i$. We induce on $X$ a binary relation $\leq$
defined by
$$
y_1\leq y_2\doteq\forall_{z\in X'}(y_1\sim z\ra y_2\sim z).
$$

In any model (even infinite) the properties (P1)--(P8) assure that $\leq$ is
a linear order with a least and greatest element. Furthermore,
every $y$ has a successor $y^+$ and a predecessor $y^-$ except when $y$
is the last or first element of $X$ respectively.

\subsubsection{Coordinatization}
We now give a formula $\COOR(x,x',t,t',z)$
that shall coordinatize the neighborhood of $z$. Let $X,X',T,T',Z$
denote the neighborhoods of $x,x',t,t',z$ respectively.
Then $\COOR$ is the conjunction of the following:
\begin{description}
\item[(C1)]
$x,x,',t,t',z,X,X',T,T',Z$ are all disjoint.
$Z$ is an independent set. All neighbors of $Z$ are in $\{z\}\cup X\cup T$.
There is no edge between $X\cup X'$ and $T\cup T'$.
\item[(C2)]
$P(x,x')\wedge P(t,t')$
\item[(C3)]
$\forall_{z\in Z}(\exists!_{x\in X}z\sim x\wedge\exists!_{t\in T} z\sim t)$
\item[(C4)]
$\forall_{x\in X}\forall_{t\in T}\exists!_{z\in Z} (z\sim x \wedge z\sim t)$
\end{description}
Thus, each $z\in Z$ has a unique pair of coordinates $(x,t)$ and
each $(x,t)$ corresponds to a unique $z$.
Note that $\qr{\COOR}=\qr P=4$ and $\alt{\COOR}=\alt P=2$.

\subsubsection{New functional and constant symbols}

To facilitate further description of $A_M$, we will use new functional
symbols. In particular, this will allow us to have new constant symbols
as symbols of nullary functions.

Writing $\bar v$, we will mean a finite sequence of variables
$v_1,v_2,\ldots$. As soon as a statement $\forall\bar y\exists!xF(x,\bar y)$
is put in $A_M$ or is derivable from what is already put in $A_M$,
we may want to denote this unique $x$ by $\phi(\bar y)$ and use
$\phi$ as a new functional symbol in the standard way. Namely,
if $Q(u,\bar z)$ is a formula with free variables $u,\bar z$, then
\begin{eqnarray*}
Q(\phi(\bar y),\bar z)&\doteq&\exists x(F(x,\bar y)\And Q(x,\bar z))
\mbox{\ \ or}\\
Q(\phi(\bar y),\bar z)&\doteq&\forall x(F(x,\bar y)\imply Q(x,\bar z)).
\end{eqnarray*}
Both variants are admissible and an appropriate choice of one of them
may reduce the alternation number of a formula.
Furthermore, in this way we can express compositions of several
functions (e.g.\ \cite[section 2.9]{Men}).

In particular, in any model of $\COOR(x,x',t,t',z)$ we let
$1,2$ denote the first two elements of $X$ (under $\leq$)
and $0$ (it will represent time zero) the first element of $T$.
The same character $\omega$ will be used for the
last element of $X$ or $T$, dependent on context.
For $v$ in $X$ or $T$, $v^-$ and $v^+$ are respectively its
predecessor and successor (when defined).
The notation $(x,t)$ will be used as a binary function symbol
with meaning as explained in the preceding subsection.

\subsection{Capturing a computation by a formula}

\subsubsection{Definition of a Turing machine}\label{ss:tmdef}

By technical reasons, we prefer to use the model of a Turing machine
where the {\em tape\/} is infinite in one direction. It is known
(e.g.\ \cite[section 41]{Kle}) that it is equivalent to the model
with the tape infinite in both directions. At the start the tape
consists of the special ``Left End of Tape'' symbol $L$,
followed by an input word written down in the binary alphabet
$\{a,b\}$, and followed onward by all ``blank'' symbols $B$.
A symbol occupies one
cell. Let $s_1,\ldots,s_k$ be {\em states\/} of a Turing machine
$M$, with $s_1$ the {\em inital\/} state and $s_k$ the {\em final\/} state.
At the start $M$ is in state $s_1$ and its head is at the first $B$.
A machine is defined by a set of instructions of the following type,
where $\ga,\gb\in\{L,a,b,B\}$.
\begin{description}
\item[$s_i\ga\gb s_j$:]
If in state $s_i$ reading a symbol $\ga$, overwrite $\gb$ and
go to state $s_j$.
\item[$s_i\ga\,\mbox{\it Right}\, s_j$:]
If in state $s_i$ reading a symbol $\ga$, move the head one cell to the right
and go to state $s_j$.
\item[$s_i\ga\,\mbox{\it Left}\, s_j$:]
If in state $s_i$ reading a symbol $\ga$, move the head one cell to the left
and go to state $s_j$.
\end{description}
If $\ga=L$ in an instruction of the first type, then $\gb=L$.
This is the only case when $\gb=L$.
There is no instruction of the third type (``move to the left'')
for $\ga=L$. With this exception, for every $i<k$ and $\ga$ there is
a unique instruction what to do in state $s_i$ reading $\ga$.
The machine halts immediately after coming to state $s_k$.
If $M$ halts, its {\em running time\/} is the number of instructions
executed before termination.

\subsubsection{Formula $A_M$}

For notation simplicity, we use the same name for variables and
corresponding semantical objects (ingredients of $M$ and vertices
of a graph $G_M$). The vertex $H$ below shall be used to keep track
of the tape header. $A_M$ is the conjunction of the two graph axioms
and a long formula of the form
$$
\exists_{x,x',t,t',z,s_1,\ldots,s_k,a,b,B,L,H}
B_M(x,x',t,t',z,s_1,\ldots,s_k,a,b,B,L,H).
$$
The formula $B_M$ whose all free variables are listed above is the
conjunction of the following subformulas, where
$X,X',T,T',Z$ denote, as before, the neighborhoods of $x,x',t,t',z$
respectively.\\
{\bf (A1)}
$x,x',t,t',z,s_1,\ldots,s_k,a,b,B,L,H,X,X',T,T',Z$ are disjoint and
consist of {\em all} the vertices of the graph.\\
{\bf (A2)}
$\COOR(x,x',t,t',z)$\\
{\bf (A3)}
All of the neighbors of $a,b,B,L,H$ are in $Z$.\\
{\bf (A4)}
For all $x\in X$ and $t\in T$ the vertex $(x,t)$ is adjacent to precisely
one of $a,b,B,L$.  We will write $\VAL(x,t)$ for this value,
which represents the symbol on the Turing Machine at position (cell of the
tape) $x$ and time $t$. Note that, as $\VAL(x,t)$ ranges over four
possible values $L,a,b,B$, using this functional symbol requites no
extra quantification. For example, the formula $\VAL(x,t)=\ga$ reads
just $(x,t)\sim\ga$.\\
{\bf (A5)}
All neighbors of $H$ are in $Z$.
For all $t\in T$ there is a unique $x\in X$ for which $(x,t)$ is adjacent
to $H$.  We write $\HP(t)$ for this $x$, which represents the header
position. Thus, $\HP(t)=x$ reads $(x,t)\sim H$.
We shall write $\VAL(t)=\VAL(\HP(t),t)$, the symbol that the header is
looking at time $t$. If $\HP$ is used within $\VAL$, it takes one extra
quantifier. Note that a subformula $\VAL(t)=\ga$ has quantifier rank 2
and alternation number 0. Furthermore, $\VAL(t^+)=\ga$ has quantifier rank
4 and can be written with alternation number 0.\\
{\bf (A6)}
The neighbors of $s_1,\ldots,s_k$ are all in $T$.  For all $t\in T$
precisely one of $s_1,\ldots,s_k$ is adjacent to $t$.
We write $\ST(t)$ for this $s_i$, which represents the state at time $t$.
Note that $\ST(t)=s\doteq t\sim s$.

We want the Turing machine to start in the standard position:\\
{\bf (A7)}
$\VAL(1,0)=L\wedge \forall_{x\neq 1}\VAL(x,0)=B \wedge \HP(0)=2
\wedge \ST(0)=s_1$\\
We want the Turing Machine to end in the final state and not be there
before that:\\
{\bf (A8)}
$\forall_{t\in T}(\ST(t)=s_k\equi t=\omega)$\\
We want values on the tape not to change except (possibly) at the header
position:\\
{\bf (A9)}
$\forall_{t\in T,t\neq \omega} \forall_{x\in X} (x\neq\HP(t)\ra
\VAL(x,t^+)=\VAL(x,t))$\\
We want the rightmost spot on the tape to be used. (We need this for
uniqueness of the model, we don't want to allow superfluous blanks.)\\
{\bf (A10)}
$\exists_{t\in T}\VAL(\omega,t)\neq B$\\
We need that the instructions would not push the Turing Machine to the right
of $x=\omega$.  For every $s_i,\ah$ such that when at state $s_i$ and
value $\ah$ the instruction push the header to the right we have:\\
{\bf (A11)}
$\neg \exists_{t\in T}(\VAL(t)=\ah \And \ST(t)=s_i \wedge HP(t)=\omega)$

We are down to the core workings of the Turing Machine.
For each instruction of the first type we have:\\
{\bf (A12)}
$\forall_{t\in T}\forall_{x\in X}(
\ST(t)=s_i\wedge \HP(t)=x\wedge \VAL(t)=\ah \ra
\ST(t^+)=s_j \wedge \VAL(t^+)=\beta \wedge \HP(t^+)=x)$\\
For each instruction of the second type we have:\\
{\bf (A13)}
$\forall_{t\in T}\forall_{x\in X}(
\ST(t)=s_i\wedge \HP(t)=x\wedge \VAL(t)=\ah \ra
\ST(t^+)=s_j \wedge \HP(t^+)=x^+ \wedge \VAL(x,t^+)=\ah)$\\
For each instruction of the third type we have:\\
{\bf (A14)}
$\forall_{t\in T}\forall_{x\in X}(
\ST(t)=s_i\wedge \HP(t)=x\wedge \VAL(t)=\ah \ra
\ST(t^+)=s_j \wedge \HP(t^+)=x^- \wedge \VAL(x,t^+)=\ah)$

\subsubsection{Proof of the Simulation Lemma}

A straightforward inspection shows that $\qr{B_M}=6$, contributed,
for example, by (A9). This gives Item 1 of the lemma.
Since we treat a variable as a single symbol, (A1) and (A6) have length
$O(k^2)$, (A11)--(A14) have length $O(k)$, and all the others have
constant length. This gives Item~2.
A straightforward inspection shows that $\alt {B_M}=2$,
contributed by (A2). This gives Item~3.

Item 4 requires a bit of extra work. As $A_M\in\altsete3$,
Lemma~\ref{lem:alttoprenex} implies that $A_M$ is reducible
to an equivalent prenex formula with quantifier prefix
$\exists^*\forall^*\exists^*\forall^*$. We make a stronger claim
that one can achieve the prefix
$\exists^*\forall^{O(1)}\exists^{O(1)}\forall^{O(1)}$. Note that
$B_M$ has a constant number of conjunctive members with constant length
and hence those contribute a constant number of quantifiers.
(A1) and (A6), though have length dependent on $k$, contain a constant
number of quantifiers. The remaining (A11)--(A14) should be tackled
with more care as every of these components, though has constant
number of quantifiers, occurs in $B_M$ in $O(k)$ variants
for various pairs $s_i,\ga$. Fortunately, all these occurrences
can be replaced by a single formula with constantly many quantifiers.
For example, introducing two new variables $s$ and $c$, we can replace
the conjunction of all variants of (A11) by
$$
\neg \exists_{t\in T}\exists_s\exists_c
[\bigvee_{s_i,\ga}(s=s_i\And c=\ga)
\And \VAL(t)=c \And \ST(t)=s \wedge HP(t)=\omega],
$$
where the disjunction is over the specified pairs $s_i,\ga$.

Let us turn to Items 5 and 6. It should be clear that, if $M$ halts,
its computation is converted to a graph satisfying $A_M$, whose order
exceeds the running time. Such a graph is unique up to isomorphism
because the adjacencies of any finite model of $A_M$ must mirror
the actions of the Turing machine. By the same reason,
any finite model of $A_M$ is converted into a halting computation
of $M$ and hence, if $A_M$
has a finite model, then $M$ halts on the empty input.
It remains to notice that, if $M$ halts, then $A_M$ has no infinite model.
Let $m$ be the running time of $M$. In any model of $A_M$, the first $m$
values of $t$ must simulate $m$ steps of $M$'s computation.
By (A8), the set $T$ is therefore finite. By (A10),
the cardinality of $X$ cannot exceed the cardinality of $T$ and hence
$X$ is finite too. It immediately follows that the other components of the
model, $X'$, $T'$, and $Z$, are finite as well. The proof is complete.

\section{Other consequences of the Simulation Lemma}\label{s:other}

\subsection{There are succinct definitions by prenex formulas}\label{ss:pre}

Due to \refeq{eq:defsimpl}, any graph of order $n$ is definable by a
prenex formula of quantifier rank $n+1$ with alternation number 1.
Though the class of prenex formulas may appear rather restrictive,
it turns out that, if one allows to increase the alternation number
to 3, then there are graphs definable by prenex formulas with very
small quantifier rank.

\begin{definition}
Let $\lpre aG$ denote the minimum length of a closed prenex formula
with alternation number at most $a$ that defines a graph $G$. Furthermore,
$$
\sfpre an = \min_{|G|=n}\lpre aG.
$$
\end{definition}

\begin{theorem}\label{thm:prenexpos}
There is no general recursive function $f$ such that
$f(\sfpre 3n)\ge n$ for all $n$.
\end{theorem}

\begin{proof}
We proceed precisely as in the proof of Theorem \ref{thm:hugegap}
but using, instead of $A_M$, the prenex formula $P_M$ given by
the Simulation Lemma. We will need a recursive bound $|P_{M_k}|\le l(k)$.
We can take $l(k)=ck^24^{k^2}$ owing to Lemma~\ref{lem:lvsdprenex}.
\end{proof}

\subsection{The set of defining sentences is undecidable}

\begin{theorem}\label{thm:defundec}
The class of defining sentences is undecidable.
\end{theorem}

\begin{proof}
Given a Turing machine $M$, consider a sentence $A_M$ as in the
Simulation Lemma. If $M$ halts on the empty input, $A_M$ is defining.
Suppose that $M$ never halts. Then either $A_M$ has no model or it
has an infinite model. By Lemma \ref{lem:findef}, $A_M$ is
not defining in both cases. We therewith have reduced the halting problem
(for the empty input) to the decision problem for the set of defining
sentences.
\end{proof}

Note a partial positive result given by Lemma~\ref{lem:defdec}.

\subsection{$D_0(G)$ and $D(G)$ are not recursively related}

Obviously, $D(G)\le D_0(G)$ for all graphs $G$.
How far apart from each other can be these two values?
Is there a converse
relation $D_0(G)\le f(D(G))$, for any general recursive function $f$?
The answer is ``no''. We will actually prove a stronger fact.
Let $D_{1/2}(G)$ denote the minimum quantifier rank of a
$\altsete1$-sentence that defines $G$. Notice the hierarchy
$$
D(G)\le D_3(G)\le D_2(G)\le D_1(G)\le D_{1/2}(G)\le D_0(G).
$$
We are able to show a superrecursive gap even between $D_3(G)$ and
$D_{1/2}(G)$.

\begin{theorem}\label{thm:dvsd0}
There is no general recursive function $f$ such that
$$
D_{1/2}(G)\le f(D_3(G))
$$
for all graphs~$G$.
\end{theorem}

\begin{lemma}\label{lem:decalt1}
The finite satisfiability of a $\altsete1$-sentence is decidable.
\end{lemma}

\begin{proof}
By Lemma \ref{lem:alttoprenex}, a $\altsete1$-sentence
effectively reduces to an equivalent formula in the \BS\/ class.
The finite satisfiability of the latter is decidable by the
Ramsey theorem.
\end{proof}

The next lemma is related to the well-known fact that, over a finite
vocabulary, there are only finitely many pairwise inequivalent
sentences of bounded quantifier rank (cf.\ \cite[lemma 4.4]{CFI}).

\begin{lemma}\label{lem:uniset}
Given $m\ge0$, one can effectively construct a finite set $U_m$
consisting of $\altsete1$-sentences of quantifier rank
$m$ so that every $\altsete1$-sentence of quantifier rank $m$
has an equivalent in~$U_m$.
\end{lemma}

\begin{proof}
Any sentence $A$ of quantifier rank $m$ can be rewritten in an equivalent
form $A'$ so that $A'$ uses at most $m$ variables, where different
occurrences of the same variable are not counted
(see e.g.\ \cite[proposition 2.3]{PVV}). Referring to this fact,
we will put in $U_m$ only sentences over the variable set
$\{x_1,\ldots,x_m\}$. We now prove the lemma in a stronger form
saying that, for each $m$ and $k$ such that $0\le k\le m$,
one can construct a finite set $U_{m,k}$ which is universal
for the class of $\altsete1$-formulas of quantifier rank $k$
over the variable set $\{x_1,\ldots,x_m\}$ with precisely $k$ variables
bound.

We proceed by induction on $k$. Consider the base case of $k=0$.
There are $a=2{m\choose 2}$ atomic formulas $x_i\sim x_j$ and $x_i=x_j$.
Any quantifier-free formula is a Boolean combination of these
and can be represented by a perfect
DNF (except the totally false formula for which we fix representation
$x_1=x_1\And x_1\ne x_1$). The set $U_{m,0}$ consists of all $2^{2^a}$
such expressions.

\newcommand{\ue}[1]{U^{\exists}_{m,#1}}
\newcommand{\uu}[1]{U^{\forall}_{m,#1}}

$U_{m,k}$ will consist of two parts, $\ue k$ and $\uu k$, the former
for formulas with at least one existential quantifier and the latter
for formulas with no existential quantifier. If $k=0$, we have
$\ue0=\emptyset$ and $\uu0=U_{m,0}$. Assume that $k\ge1$ and $U_{m,k-1}$
has been already constructed. We construct $U_{m,k}$ in four steps.
\begin{enumerate}
\item
Put in $\ue k$ the formulas $\exists x_i A$ for all
$A\in U_{m,k-1}$ and $i\le m$ such that no occurrence of $x_i$ in $A$
is bound.
\item
Put in $\uu k$ the formulas $\forall x_i A$ for all
$A\in \uu{k-1}$ and $i\le m$ such that no occurrence of $x_i$ in $A$
is bound.
\item
Put in $\ue k$ all monotone Boolean combinations of formulas
from $\ue k$ and $\uu k$ as constructed in Steps 1 and 2
with at least one formula from $\ue k$ involved.
\item
Put in $\uu k$ all monotone Boolean combinations of formulas
from $\uu k$ as constructed in Step~2.
\end{enumerate}
Finally, to obtain $U_m$ exactly as claimed in the lemma, we
set $U_m=U_{m,m}$.
\end{proof}

\begin{proofof}{Theorem \ref{thm:dvsd0}}
Suppose on the contrary that a such $f$ exists. Using the $f$,
we will design an algorithm for the halting problem, contradicting
the unsolvability of the latter.

Given a Turing machine $M$, we construct the sentence $A_M$ as in
the Simulation Lemma. Recall that
\begin{itemize}
\item
$\alt {A_M}=3$;
\item
if $M$ halts on the empty input, then $A_M$ defines a finite graph $G_M$;
\item
if $M$ does not halt, then $A_M$ has no finite model.
\end{itemize}
Denote $k=\qr {A_M}$ and $m=\max_{i\le k}f(i)$. Thus, if $G_M$ exists,
then $D_3(G_M)\le k$ and, by the assumption, $D_{1/2}(G_M)\le m$.

Construct $U_m$ as in Lemma \ref{lem:uniset}
and add to every sentence in $U_m$ the two graph axioms. We know that
$U_m$ contains a sentence defining $G_M$ and this will help us to
construct this graph (if it exists). Remove from $U_m$ all finitely
unsatisfiable formulas. This task is tractable by Lemma \ref{lem:decalt1}.
For every remaining sentence, by brute-force search we eventually find
a finite graph satisfying it (we need one model for every sentence
and do not care that some sentences may have other models).
Let $G_1,\ldots,G_l$ be the list of these graphs.

If $M$ halts, one of the $G_i$'s coincides with $G_M$ and satisfies $A_M$.
If $M$ does not, none of the $G_i$'s satisfies $A_M$. Thus, the verification
if $A_M$ is true on one of the $G_i$'s allows us to recognize if
$M$ halts on the empty input.
\end{proofof}

\begin{corollary}\label{cor:dvsd0}
\mbox{}

\begin{enumerate}
\item
There is no general recursive function $f$ such that
$D_0(G)\le f(D(G))$ for all graphs~$G$.
\item
There is no general recursive function $f$ such that
$D_0(G,G')\le f(D(G,G'))$ for all non-isomorphic $G$ and~$G'$.
\end{enumerate}
\end{corollary}

\begin{proof}
1) Suppose on the contrary that a such $f$ exists. Then we would have
$D_{1/2}(G)\le D_0(G)\le f(D(G))\le\max_{i\le D_3(G)}f(i)$,
contradictory to Theorem~\ref{thm:dvsd0}.

2) Again, suppose that a such $f$ exists. By Proposition \ref{prop:ddd},
$D_0(G)=D_0(G,G')$ for some $G'$. It follows that
$D_0(G)\le f(D(G,G'))\le\max_{i\le D(G)}f(i)$, contradictory to Item~1.
\end{proof}

It is also worthy to note the following fact.

\begin{theorem}\label{thm:comput0}
$D_0(G)$ and $D_{1/2}(G)$ are computable functions of graphs.
\end{theorem}

\begin{proof}
We prove the theorem for $D_{1/2}(G)$; For $D_0(G)$ the proof is similar.
Starting from $m=2$, we trace through the universal set $U_m$ given
by Lemma \ref{lem:uniset} and, for each sentence $A\in U_m$,
check whether $G$ satisfies $A$ and, if so, whether $A$ is defining.
The latter can be done on the account of Lemma \ref{lem:defdec}.
If no such $A$ is found, we conclude that $D_{1/2}(G)>m$ and increase $m$
by~1.
\end{proof}

\begin{remark}
A variant of Theorem \ref{thm:dvsd0} for the
formula length is also true, even with a simpler proof
(no reference to Lemma \ref{lem:uniset} is needed).
\end{remark}

\subsection{An undecidable fragment of the theory of finite graphs}%
\label{ss:tfg}

Given a class of $\vocab$-structures $\calC$, let $\sat\calC$
(resp.\ $\sateq\calC$) be the set of formulas over $\vocab$
without equality (resp.\ with equality) that have a model in $\calC$.
Furthermore, let $\finsat\calC$ (resp.\ $\finsateq\calC$) be the set of
formulas over $\vocab$ without equality (resp.\ with equality) that
have a finite model in $\calC$.
If $X$ is one of the aforementioned sets and $F$ is a class of
formulas over $\vocab$, we call the intersection $F\cap X$
the {\em $F$-fragment of\/} $X$. We will be interested in the case
that $F$ is a {\em prefix class}, that is, consists of prenex formulas
whose quantifier prefix agrees with a given pattern. Describing such
a pattern, we use $\forall^*$ or $\exists^*$ to denote a
string of all $\forall$ or all $\exists$ of any length.

Let $\calD$ (resp.\ $\calI$ and $\calS$) denote the class of structures
consisting of a single binary relation (resp.\ irreflexive binary relation
and symmetric binary relation). In other words, $\calD$ is the class
of directed graphs.
By $\calG$ we denote the class of graphs, i.e., structures consisting
of a single irreflexive symmetric relation.

Based on Church's and Turing's solution of Hilbert's
{\em Entscheidungsproblem}, Kalm\'ar \cite{Kal} proved that
$\sat\calD$ is undecidable. Following the Kalm\'ar result and
the Trakhtenbrot theorem \cite{Tra}, Vaught \cite{Vau} proved that
the set $\finsat\calD$ and the set of formulas not in $\sat\calD$
are recursively inseparable, that is, no decidable
set contains the former and is disjoint with the latter.
In particular, both $\finsat\calD$ and $\sat\calD$ are undecidable.
Currently a complete classification of prefix fragments of
$\sat\calD$, $\finsat\calD$, $\sateq\calD$, and $\finsateq\calD$
is known (see \cite{BGG}, a reference book on the subject).

Gurevich \cite{Gur} shows that the $\forall^3\exists^*$-fragment of
$\sat\calI$ is undecidable while the $\forall^2\exists^*$-fragment
is decidable (cf.\ \cite[corollary 6.2.35]{BGG}).
Church and Quine \cite{CQu} established the undecidability of
$\sat\calS$. Note that this result is easily extended to
$\sateq\calG$. The undecidability of $\sat\calG$ was proved by
Rogers \cite{Rog}. Lavrov \cite{Lav} (see also \cite[theorem 3.3.3]{ELTT})
improved this by showing the recursive inseparability of $\finsat\calG$ and
the set of formulas not in $\sat\calG$.

Lavrov's proof provides us with a reduction of the decision problem for
$\calD$ to the decision problem for $\calG$. If combined with
the known results on undecidable fragments of $\finsat\calD$, this gives
us some undecidable fragments of $\finsat\calG$, for example,
$\forall^9\exists^*\forall^*\exists^*$. However, this method
apparently cannot give undecidable fragments with less than two
star symbols.
Our Simulation Lemma has relevance to this circle
of questions.

\begin{theorem}\label{thm:fgtundec}
For some $l$, $m$, and $n$, the
$\exists^*\forall^l\exists^m\forall^n$-fragment of $\finsateq\calG$
is undecidable.
\end{theorem}

\begin{proof}
By the Simulation Lemma, a Turing machine $M$ halts on the empty input
iff the formula $A_M$ has a finite graph as a model.
Thus, the conversion of $A_M$
to a prenex formula according to Item 4 of the Simulation Lemma reduces
this variant of the halting problem to the satisfiability problem
for $\exists^*\forall^l\exists^m\forall^n$-formulas over finite graphs.
\end{proof}

The theorem should be contrasted with the decidability of the
$\exists^*\forall^*$-fragment, which follows from the Ramsey theorem
and the fact that the class of graphs is definable by a
$\exists^*\forall^*$-formula.
We do not try to specify numbers $l,m,n$ since the values derivable
from our proof are, though not so big, surely improvable by extra technical
efforts. Note that a variant of the theorem for $\finsat\calD$
is known to be true with best possible $l=m=n=1$
(see \cite[theorem 3.3.2]{BGG}, which is Sur\'anyi's theorem extended
to the finite satisfiability by Gurevich).

Note another equivalent form of Theorem \ref{thm:fgtundec}.
Let $\fintheq\calG$ denote the {\em first order theory of finite graphs
with equality}, i.e., the set of first order sentences with relation
symbols $\sim$ and $=$ that are true on all finite graphs.
Observe that a sentence $A$ is in $\fintheq\calG$ iff $\neg A$ is not
in $\finsateq\calG$. It follows that the
$\forall^*\exists^l\forall^m\exists^n$-fragment
of $\fintheq\calG$ is undecidable.

\section{The succinctness function over trees: Upper bound}\label{s:tr_up}

We define a variant of the succinctness function for a class of
graphs $\cal C$ (with respect to the quantifier rank) by
$$
\sfclass=\min\setdef{D(G)}{G\in{\cal C},\,|G|=n}.
$$
We here prove a log-star upper bound for the class of trees.

\begin{theorem}\label{thm:tr_up}
$\sftree<\log^*n+5$.
\end{theorem}

\noindent
The proof takes the rest of this section.

\subsection{Rooted trees}

A {\em rooted tree\/} is a tree with one distinguished vertex, which is
called the {\em root}. If $T$ is a tree and $v\in V(T)$, then $T_v$
denotes the tree $T$ rooted at $v$. An isomorphism of rooted trees
should not only preserve the adjacency relation but also map one root
to the other. Thus, for distinct $u,v\in V(T)$, rooted trees $T_u$
and $T_v$, though have the same underlying tree $T$, may be non-isomorphic.

An {\em automorphism\/} of a rooted tree is an isomorphism
from it onto itself.
Obviously, any automorphism leaves the root fixed. We call a rooted tree
{\em asymmetric\/} if it has no non-trivial automorphisms, that is,
no automorphisms except the identity.

The {\em depth\/} of a rooted tree $T_v$, which is denoted by $\depth T_v$,
is the eccentricity of its root.
If $(v,\ldots,u,w)$ is a path in $T_v$, then $w$ is called a {\em child\/}
of $u$.
We define the relation of being a {\em descendant\/} to be the
transitive and reflexive closure of the relation of being a child.

If $w\in V(T_v)$, then $T_v(w)$ denotes the subtree of $T_v$ spanned
by the set of all descendants of $w$ and rooted at $w$.
If $w$ is a child of $u\in V(T_v)$, then $T_v(w)$ is called
a {\em $u$-branch\/} of $T_v$.

\subsection{Diverging trees}

We call $T_v$ {\em diverging\/} if, for every vertex $u\in V(T_v)$,
all $u$-branches of $T_v$ are pairwise non-isomorphic.

\begin{lemma}\label{lem:down}
A rooted tree $T_v$ is diverging iff its $v$-branches are pairwise
non-isomorphic and each of them is diverging.
\end{lemma}

\begin{proof}
Assume that $T_v$ is diverging. Its $v$-branches are pairwise non-isomorphic
by the definition. Furthermore, let $T_v(w)$ be a $v$-branch of $T_v$ and
$u\in V(T_v(w))$. Note that any $u$-branch of $T_v(w)$ is also a
$u$-branch of $T_v$. Therefore, all of them are pairwise
non-isomorphic and $T_v(w)$ is diverging.

For the other direction, consider a non-root vertex $u$ of $T_v$
and let $T_v(w)$ be the $v$-branch of $T_v$ containing $u$ ($w=u$ is possible).
Note that any $u$-branch of $T_v$ is also a
$u$-branch of $T_v(w)$. Therefore, all of them are pairwise
non-isomorphic and we conclude that $T_v$ is diverging.
\end{proof}

\begin{lemma}\label{lem:divvsasym}
A rooted tree $T_v$ is diverging iff it is asymmetric.
\end{lemma}

\begin{proof}
We proceed by induction on $d=\depth T_v$. The base case of $d=0$
is trivial. Let $d\ge1$.

Assume that $T_v$ is diverging. By Lemma \ref{lem:down},
no automorphism of $T_v$ can map one $v$-branch onto another $v$-branch.
By the same lemma and the induction assumption, no non-trivial
automorphism can map a $v$-branch onto itself. Thus, $T_v$ has no
non-trivial automorphism.

Assume now that $T_v$ is asymmetric. Hence all $v$-branches are
pairwise non-isomorphic and each of them is asymmetric.
By the induction assumption, each $v$-branch is diverging.
By Lemma \ref{lem:down} we conclude that $T_v$ is diverging.
\end{proof}

We now carry over the notion of a diverging tree to (unrooted) trees.
Clearly, any automorphism of a tree $T$ either leaves central vertices
$c_1$ and $c_2$ fixed or transposes them ($c_1=c_2$ if the diameter
$d(T)$ is even). If $d(T)$ is odd, Lemma \ref{lem:divvsasym}
implies that $T_{c_1}$ and $T_{c_2}$ are simultaneously diverging or not.
This makes the following definition correct: A tree $T$ is {\em diverging\/}
if the rooted tree $T_c$ for a central vertex $c$ is diverging.
It is not hard to see that $T$ is diverging iff one of the following
conditions is met:
\begin{enumerate}
\item
$T$ has no non-trivial automorphism.
\item
$T$ has exactly one non-trivial automorphism and this automorphism
transposes two central vertices of~$T$.
\end{enumerate}

\subsection{Spoiler's strategy}

In this section we exploit the characterization of the quantifier rank
of a distinguishing formula as the length of the \EF\/ game
(see Proposition~\ref{prop:loggames}).

\begin{lemma}\label{lem:distance}
Suppose that in the \EF\/ game on $(G,G')$ some two
vertices $x,y\in V(G)$ at distance $k$ were selected so that their
counterparts $x',y'\in V(G')$ are at a strictly larger distance
(possibly infinity).

Then Spoiler can win in at most $\lceil\log k\rceil$ extra moves, playing
all the time inside~$G$.
\end{lemma}

\begin{proof}
Spoiler sets $u_1=x$, $u_2=y$, $v_1=x'$, $v_2=y'$, and places
a pebble on the middle vertex $u$ in a shortest path from $u_1$ to $u_2$
(or either of the two middle vertices if $d(u_1,u_2)$ is odd).
Let $v\in V(G')$ be selected by Duplicator in response to $u$.
By the triangle inequality, we have $d(u,u_m)<d(v,v_m)$
for $m=1$ or $m=2$. For such $m$ Spoiler resets $u_1=u$, $u_2=u_m$,
$v_1=v$, $v_2=v_m$ and applies
the same strategy once again. Therewith Spoiler ensures that, in each round,
$d(u_1,u_2)<d(v_1,v_2)$. Eventually, unless Duplicator loses earlier,
$d(u_1,u_2)=1$ while $d(v_1,v_2)>1$, that is, Duplicator fails to
preserve adjacency.

To estimate the number of moves made, notice that initially
$d(u_1,u_2)=k$ and for each subsequent $u_1,u_2$ this distance
becomes at most $f(d(u_1,u_2))$, where $f(\ga)=(\ga+1)/2$.
Therefore the number of moves does not exceed the minimum $i$
such that $f^{(i)}(k)<2$. As $(f^{(i)})^{-1}(\gb)=2^i\gb-2^i+1$,
the latter inequality is equivalent to
$2^i\ge k$, which proves the bound.
\end{proof}

Note that the bound of Lemma \ref{lem:distance} is tight,
more precisely, it cannot be improved to $\lceil\log k\rceil-1$.
For example, let $C_n$ denote a cycle of length $n$ and
$2C_n$ the disjoint union of two such cycles.
It is known
(e.g.\ \cite[Proof of Theorem 2.4.2]{Spe} or \cite[Example 2.3.8]{EFl})
that Duplicator can survive in the \EF\/ game on $C_{2k+1}$ and $C_{2k+2}$
in more than $\log k+1$ rounds for any strategy of Spoiler,
in particular, when Spoiler begins with selecting two antipodal
vertices in $C_{2k+2}$.
Furthermore, if $d(x',y')=\infty$, Duplicator can be persistent as well.
For example, she can survive in the game on $C_{2k}$ and $2C_{2k}$
during $\lfloor\log(2k-1)\rfloor$ rounds for any strategy of Spoiler,
in particular, when Spoiler's first move is in one component of $2C_{2k}$
and his second move is in the other component of $2C_{2k}$
(e.g.\ \cite[Example 2.3.8]{EFl}).

\begin{lemma}\label{lem:distdiam}
If graphs $G$ and $G'$ have different diameters
(including the case that $G$ is connected and $G'$ is disconnected),
then $D_1(G,G')\le\lceil\log d(G)\rceil+2$.
\end{lemma}

\begin{proof}
Assume that $d(G)<d(G')$.
Spoiler begins with selecting two vertices at distance $d(G)+1$ in $G'$,
then jumps to $G$, and uses the strategy of Lemma~\ref{lem:distance}.
\end{proof}

\begin{lemma}\label{lem:nontree}
If $G$ is a tree, $G'$ is a connected non-tree, and $d(G)=d(G')$, then
$D_0(G,G')<\lceil\log d(G)\rceil+4$.
\end{lemma}

\begin{proof}
Denote $k=d(G)=d(G')$.
Let $C$ be a shortest cycle in $G'$. Notice that $C$ has length
at most $2k+1$. Spoiler begins with selecting in $C$ a vertex $z'$
along with its neighbors $x'$ and $y'$.
Let $z$, $x$, and $y$ be the corresponding responses of Duplicator in $G$.
The vertex $z$ cannot be a leaf of $G$ for else Duplicator has lost.
{}From now on Spoiler plays all the time in $H'=G'-z'$
and Duplicator is enforced to play in $H=G-z$.
In these graphs $d(x',y')\le 2k-1$ and $d(x,y)=\infty$.
Therefore the strategy of Lemma \ref{lem:distance} applies and
Spoiler wins in at most $\lceil\log(2k-1)\rceil$ extra moves.
\end{proof}

\begin{lemma}\label{lem:bothdiv}
Let $T$ and $T'$ be two non-isomorphic diverging trees with $d(T)=d(T')$
(and hence $r(T)=r(T')$). Then $D(T,T')\le r(T)+1$.
\end{lemma}

\begin{proof}
In the first move Spoiler selects $x$, a central vertex of $T$.
Duplicator's response, $x'$, should be a central vertex of $T'$
because otherwise Spoiler selects a vertex $y'$ in $T'$ with
$d(x',y')>r(T)$ and applies the strategy of Lemma \ref{lem:distance}.
We will denote the vertices selected by the players in $T$ and $T'$
during the $i$-th round by $x_i$ and $x'_i$; In particular,
$x_1=x$ and $x'_1=x'$. Spoiler will play so
that $(x_1,\ldots,x_i)$ and $(x'_1,\ldots,x'_i)$ are always paths.
Another condition that will be obeyed by Spoiler
is that $T_x(x_i)$ and $T'_{x'}(x'_i)$ are non-isomorphic.

Assume that the $i$-th round has been played. If exactly one of
the vertices $x_i$ and $x'_i$ is a leaf (we will call a such
situation terminal), then Spoiler prolongs that path for which this
is possible and wins. Assume that neither of $x_i$ and $x'_i$ is a leaf
and that $T_x(x_i)$ and $T'_{x'}(x'_i)$ are non-isomorphic
(in particular, this is so for $i=1$). By the definition of
a diverging rooted tree, all $T_x(u)$ with $u$ a child
of $x_i$ are pairwise non-isomorphic. The same concerns all $T'_{x'}(u')$
with $u'$ a child of $x'_i$.
It follows that there is a $T_x(u)$ not isomorphic
to any of the $T'_{x'}(u')$'s or there is a $T'_{x'}(u')$ not
isomorphic to any of the $T_x(u)$'s. Spoiler selects such $u$ for $x_{i+1}$
or $u'$ for $x'_{i+1}$. Clearly, Spoiler has an appropriate move
until a terminal situation occurs. The latter occurs in the $r(T)$-th
round at latest.
\end{proof}

\begin{lemma}\label{lem:divandnot}
Let $T$ and $T'$ be two trees with $d(T)=d(T')$
(and hence $r(T)=r(T')$). Suppose that $T$ is diverging but $T'$
is not. Then $D(T,T')\le r(T)+2$.
\end{lemma}

\begin{proof}
In the first move Spoiler selects $x'$, a central vertex of $T'$.
Similarly to the preceding proof, we may suppose that
Duplicator's response $x$ is a central vertex of $T$.
Let $y'$ be a vertex of $T'$ such that $T'_{x'}(y')$ is not diverging
but, for any child $z'$ of $y'$, $T'_{x'}(z')$ is.
Note that $y'$ must have two children $z'_1$ and $z'_2$
such that $T'_{x'}(z_1')$ and $T'_{x'}(z_2')$ are isomorphic.

In subsequent moves Spoiler selects the path $P'=(x',\ldots,y',z'_1)$.
Let $P=(x,\ldots,y,z)$ be Duplicator's response in $T$. If $T_x(z)$
and $T_{x'}(z'_1)$ have different depths $d$ and $d'$, say $d>d'$,
then Spoiler prolongs $P$ with $d'+1$ new vertices and wins.
It is clear that the prolonged path has at most $r(T)+1$ vertices.

Suppose now that $d=d'$. If $T_x(z)$ and $T_{x'}(z'_1)$ are non-isomorphic,
then Spoiler adopts the strategy of Lemma \ref{lem:bothdiv} and wins
having made totally at most $r(T)+1$ moves.
If $T_{x}(z)$ and $T'_{x'}(z_1')$ are isomorphic, then Spoiler
selects $z'_2$. In response Duplicator must select a child
of $y$ different from $z$. Denote it by $z^*$. The subtree $T_x(z^*)$
is non-isomorphic to $T_{x}(z)$ and hence to $T'_{x'}(z_2')$.
Now Spoiler is able to proceed with $T_x(z^*)$ and $T'_{x'}(z_2')$
as it was described and wins having made totally at most $r(T)+2$ moves
(one extra move was made to switch from $z'_1$ to $z'_2$).
\end{proof}

\begin{lemma}\label{lem:defdiv}
Let $T$ be a diverging tree of radius at least 6. Then $D(T)\le r(T)+2$.
\end{lemma}

\begin{proof}
Let $T'$ be a graph non-isomorphic to $T$. The pair $T,T'$ satisfies
the condition of one of Lemmas \ref{lem:distdiam}--%
\ref{lem:divandnot}. These lemmas provide us
with bound $D(T,T')\le r(T)+2$. By Proposition \ref{prop:ddd}, we
therewith have the bound for $D(T)$.
\end{proof}

We have shown that diverging trees are definable with quantifier rank
no much larger than the radius. It remains to show that, given the radius,
there are diverging trees with large order and, moreover, the orders
of these large trees fill long segments of integers.

\begin{lemma}\label{lem:number}
Given $i\ge0$, let $M_i$ denote the total number of
(pairwise non-isomorphic) diverging rooted trees of depth at most~$i$.
Then $M_i=T(i)$.
\end{lemma}

\begin{proof}
Let $m_i$ denote the number of diverging rooted trees of depth precisely
$i$. Thus, $m_0=1$ and $M_i=m_0+\ldots+m_i$.
By Lemma \ref{lem:down}, a depth-$(i+1)$
tree $T_v$ is uniquely determined by the set of its $v$-branches,
which are diverging rooted trees of depth at most $i$. Vise versa,
any set of diverging rooted trees of depth at most $i$ with at least
one tree of depth precisely $i$, determines a depth-$(i+1)$ tree.
It follows that $m_{i+1}=(2^{m_i}-1)2^{M_{i-1}}$, where we put $M_{-1}=0$.
By induction, we obtain $m_i=T(i)-T(i-1)$ and $M_i=T(i)$.
\end{proof}

Note that a diverging rooted tree of depth $i$ can have the minimum
possible number of vertices $i+1$ (a path).

\begin{lemma}\label{lem:order}
Let $N_i$ denote the maximum order of a diverging rooted tree of depth
$i$. Then $N_i>T(i-1)$.
\end{lemma}

\begin{proof}
The largest diverging rooted tree $T_v$ of depth $i$ has every
of $M_{i-1}$ diverging rooted trees of depth at most $i-1$
as a $v$-branch. Thus, $N_i>M_{i-1}=T(i-1)$.
\end{proof}

\begin{lemma}\label{lem:dense}
For every $n$ such that $i+1\le n\le N_i$ there is a diverging
rooted tree of depth $i$ and order $n$.
\end{lemma}

\begin{proof}
We proceed by induction on $i$. The base case of $i=0$ is trivial.
Let $i\ge 1$. For $n=i+1$ we are done with a path. We will prove that
any diverging rooted tree $T_v$ of depth $i$ except the path can be modified
so that it remains a diverging rooted tree of the same depth but the order
becomes 1 smaller.

Let $l$ be the smallest depth of a $v$-branch of $T_v$ and fix a branch
$T_v(w)$ of this depth with minimal order. If $T_v(w)$ is a path, we
delete its leaf. If not, we reduce it by the induction assumption.
\end{proof}

\begin{lemma}\label{lem:manytrees}
Let $i\ge2$.
For every $n$ such that $2i+2\le n\le2N_i$, there is a diverging tree
of order $n$ and radius $i+1$.
\end{lemma}

\begin{proof}
If $n=2m$ is even, consider the diverging rooted tree $T_c$ with two
$c$-branches, one of order $m$, the other of order $m-1$, and both of
depth $i$ (excepting the case that $n=2i+2$ when the smaller branch has
depth $i-1$).
Such branches do exists by Lemma \ref{lem:dense}.
If $n=2m+1$ is odd, we add the third single-vertex $c$-branch.
Since the root $c$ is a central vertex of the underlying tree,
the latter is diverging.
\end{proof}

\begin{proofof}{Theorem \ref{thm:tr_up}}
Let $n>32=2T(3)$ and let $i\ge3$ be such that $2T(i)<n\le2T(i+1)$.
By Lemma \ref{lem:order}, we have $2i+6<n<2N_{i+2}$.
Owing to Lemma \ref{lem:manytrees},
there exists a diverging tree $T$ of order $n$ and radius $i+3$.
Lemma \ref{lem:defdiv} gives $D(T)\le i+5<\log^*n+5$.

For every $n\le 32$ the required bound is provided by $P_n$, the path
on $n$ vertices. It is not hard to derive from Lemma \ref{lem:distdiam}
that $D_1(P_n)<\log n+3$ for all $n$,
which satisfies our needs for $n$ in the range.
\end{proofof}

\section{The succinctness function over trees:
Zero alternations}\label{s:alt0}

Theorem \ref{thm:tr_up} assumes no restriction on the alternation number.
We now prove a weaker analog of this theorem for
$\sftreez=\min_{|T|=n}D_0(T)$, the succinctness function over trees
with the strongest restriction on the alternation number.
This is somewhat surprising in view of
Corollary \ref{cor:dvsd0} (1) asserting that $D_0(G)$ and $D(G)$ may be
very far apart from one another.

\begin{theorem}\label{thm:alt0}
For infinitely many $n$ we have $\sftreez\le2\log^*n+O(1)$.
\end{theorem}

\noindent
The proof takes the rest of the section.

\subsection{Ranked trees}

We will modify the approach worked out in the preceding section.
The proof of Theorem \ref{thm:tr_up} was based on
Lemmas \ref{lem:distdiam}--\ref{lem:divandnot}. Note that the
alternation number in Lemma \ref{lem:nontree} is 0. In Lemma
\ref{lem:distdiam} it is 1, but the bound of this lemma is actually
stronger than we need and, at the cost of some relaxation, we will be able
to improve the alternation number to 0
(see Lemmas \ref{lem:conndistdiam} and \ref{lem:disconn} below).
The real source of non-constant alternation number is Lemma \ref{lem:bothdiv}
(Lemma \ref{lem:divandnot} reduces to Lemma \ref{lem:bothdiv} and
itself makes no new complication).
To tackle the problem, we restrict the class of diverging trees so that
we will still have relation $D_0(T)=O(r(T))$ and there will still
exist trees with $\mbox{\it Tower\,}(r(T)-O(1))$ vertices.

We begin with introducing some notions and notation concerning rooted trees.
Given a rooted tree $T_v$, let $B(T_v)$ denote the set of all $v$-branches
of $T_v$. Given rooted trees $T_1,\ldots,T_m$, we
define $T=T_1\odot\cdots\odot T_m$ to be the rooted tree
with $B(T)=\{T_1,\ldots,T_m\}$. By Lemma \ref{lem:down},
if all $T_i$ are pairwise non-isomorphic and diverging, then $T$
is diverging as well. Obviously, $\depth T=1+\max_i\depth T_i$.

Let $T'_{v'}$ and $T_v$ be rooted trees. We call $T'_{v'}$ a {\em rooted
subtree\/} of $T_v$ if $v'=v$ and $V(T')\subseteq V(T)$.

\begin{figure}
\centerline{
\unitlength=1.00mm
\special{em:linewidth 0.4pt}
\linethickness{0.4pt}
\begin{picture}(106.00,46.00)
\emline{20.00}{5.00}{1}{10.00}{45.00}{2}
\emline{10.00}{45.00}{3}{10.00}{45.00}{4}
\emline{20.00}{5.00}{5}{30.00}{35.00}{6}
\put(20.00,5.00){\circle*{2.00}}
\put(18.00,14.00){\circle*{2.00}}
\put(15.00,24.00){\circle*{2.00}}
\put(13.00,34.00){\circle*{2.00}}
\put(10.00,45.00){\circle*{2.00}}
\put(30.00,35.00){\circle*{2.00}}
\put(26.00,24.00){\circle*{2.00}}
\put(23.00,14.00){\circle*{2.00}}
\emline{40.00}{45.00}{7}{50.00}{5.00}{8}
\emline{47.00}{16.00}{9}{54.00}{24.00}{10}
\put(50.00,5.00){\circle*{2.00}}
\put(47.00,16.00){\circle*{2.00}}
\put(54.00,24.00){\circle*{2.00}}
\put(45.00,24.00){\circle*{2.00}}
\put(42.00,35.00){\circle*{2.00}}
\put(40.00,45.00){\circle*{2.00}}
\emline{65.00}{45.00}{11}{75.00}{5.00}{12}
\emline{70.00}{24.00}{13}{77.00}{35.00}{14}
\put(75.00,5.00){\circle*{2.00}}
\put(72.00,15.00){\circle*{2.00}}
\put(70.00,24.00){\circle*{2.00}}
\put(77.00,35.00){\circle*{2.00}}
\put(68.00,34.00){\circle*{2.00}}
\put(65.00,45.00){\circle*{2.00}}
\emline{85.00}{45.00}{15}{100.00}{5.00}{16}
\emline{100.00}{5.00}{17}{100.00}{25.00}{18}
\emline{100.00}{5.00}{19}{105.00}{15.00}{20}
\put(100.00,5.00){\circle*{2.00}}
\put(105.00,15.00){\circle*{2.00}}
\put(100.00,15.00){\circle*{2.00}}
\put(96.00,15.00){\circle*{2.00}}
\put(100.00,25.00){\circle*{2.00}}
\put(93.00,24.00){\circle*{2.00}}
\put(89.00,35.00){\circle*{2.00}}
\put(85.00,45.00){\circle*{2.00}}
\end{picture}
}
\caption{$R^*_0$.}
\end{figure}

For each $i\ge0$, we now define the class of rooted trees $R^*_i$ as follows.
Let $R_0^*=\{T^*_1,T^*_2,T^*_3,T^*_4\}$, the set of four rooted trees depicted
in Figure 1. Observe the following properties of this set.
\begin{description}
\item[(Z1)]
$|T^*_i|\le8$ for all $i$.
\item[(Z2)]
$\depth T^*_i=4$ for all $i$.
\item[(Z3)]
All $T^*_i$ are diverging.
\item[(Z4)]
No $T^*_i$ is isomorphic to a rooted subtree of any other~$T^*_j$.
\end{description}

Assume that $R^*_{i-1}$ is already specified.
We will need a large enough $F_i\subset 2^{R^*_{i-1}}$,
a family of subsets of $R^*_{i-1}$ which is an antichain with
respect to the inclusion (i.e.\ no member of $F_i$ is included in
any other member of $F_i$). As one of suitable possibilities
(which actually maximizes $|F_i|$ by Sperner's theorem), we fix
$$
F_i={R^*_{i-1}\choose\lfloor |R^*_{i-1}|/2\rfloor},
$$
the family of all $\lfloor |R^*_{i-1}|/2\rfloor$-element subsets of
$R^*_{i-1}$. Now
$$
R^*_i=\setdef{\bigodot_{T\in S}T}{S\in F_i}.
$$
Note that $|R^*_i|=|F_i|$.

It is clear that, if $T\in R^*_i$, then $B(T)$ consists of pairwise
non-isomorphic rooted trees in $R^*_{i-1}$. By easy induction, we
have the following properties of the class $R^*_i$ for $i\ge1$.
\begin{description}
\item[(R1)]
If $T\in R^*_i$, then $r(T)=\depth T=i+4$.
\item[(R2)]
If $T\in R^*_i$, then $d(T)=2i+8$.
\item[(R3)]
If $T\in R^*_i$, then the central vertex of $T$ is equal to the root.
\item[(R4)]
All $T\in R^*_i$ are diverging.
\item[(R5)]
If $T$ and $T'$ are different members of $R^*_i$, then neither
$B(T)\subset B(T')$ nor $B(T')\subset B(T)$.
\end{description}

We define $R_i$ to be the set of underlying trees of rooted trees in $R^*_i$.
Note that for different $T,T'\in R^*_i$ their underlying trees are
non-isomorphic. If $i=0$, this is evident. If $i\ge1$, we use the fact
that, as any isomorphism between the unrooted trees takes one central
vertex to the other, it is also an isomorphism between the rooted trees.
Note also that trees in $R_i$ are diverging.

We will call trees in $R=\bigcup_{i=1}^\infty R_i$ {\em ranked}.
If $T\in R_i$, we will say that $T$ has {\em rank $i$\/} and write
$\rk T=i$.

\begin{lemma}\label{lem:minorder}
Let $N_i$ denote the minimum order of a tree of rank $i$. Then
$N_i\ge T(i-O(1))$.
\end{lemma}

\begin{proof}
Denote $M_i=|R_i|$. By the construction, we have
$$
M_0=4,\ \ M_{i+1}={M_i\choose\lfloor M_i/2\rfloor}=
\sqrt{\frac{2+o(1)}{\pi M_i}}\,2^{M_i},
$$
and
$$
N_{i+1}>1+{M_i\choose\lfloor M_i/2\rfloor}>M_{i+1}.
$$
The lemma follows by simple estimation.
\end{proof}

\subsection{Spoiler's strategy}

Consider the \EF\/ game on rooted trees $(T_v,T'_{v'})$.
Let $x_i$ denote the vertex of $T_v$ selected in the $i$-th round.
We call a strategy for Spoiler {\em continuous\/} if he plays all the time
in $T_v$ and, for each $i$, the induced subgraph $T[\{v,x_1,\ldots,x_i\}]$
is connected.

\begin{lemma}\label{lem:bothri*}
Let $T_v$ and $T'_{v'}$ be non-isomorphic rooted trees in $R^*_i$.
Then Spoiler has a continuous winning strategy in
$\game_{i+7}(T_v,T'_{v'})$ and hence $D_0(T_v,T'_{v'})\le i+7$.
\end{lemma}

\begin{proof}
We proceed by induction on $i$. In the base case of $i=0$,
Spoiler selects all non-root vertices of $T_v$ in a continuous manner
and wins by Property (Z4). Let $i\ge1$. In the first move Spoiler
selects $w$, a child of $v$ such that, for any $w'$, child of $v'$,
branches $T_v(w)$ and $T'_{v'}(w')$ are not isomorphic. This is possible
owing to Property (R5). Let $w'$ denote Duplicator's response. Both
$T_v(w)$ and $T'_{v'}(w')$ have rank $i-1$. Spoiler now invokes
a continuous strategy winning $\game_{i+6}(T_v(w),T'_{v'}(w'))$,
which exists by the induction assumption.
\end{proof}

\begin{lemma}\label{lem:new}
Let $T$, $T'$ be trees of the same even diameter and $v$, $v'$ be their
central vertices. Assume that Spoiler selects $v$ but Duplicator
responds with a vertex different from $v'$. Then Spoiler is able to win
in the next $d(T)$ moves, playing all the time in~$T$.
\end{lemma}

\begin{proof}
In a continuous manner, Spoiler selects the vertices of a diametral
path in $T$. Let $u\ne v'$ be the vertex selected by Duplicator in
response to $v$. Duplicator should now to exhibit a path of length
$d(T')=d(T)$ with $u$ at the middle, which is impossible by
Proposition~\ref{prop:ore}.
\end{proof}

\begin{lemma}\label{lem:bothri}
Let $T$ and $T'$ be non-isomorphic ranked trees of the same rank.
Then $D_0(T,T')\le2\rk T+9$.
\end{lemma}

\begin{proof}
Let $v$ and $v'$ be central vertices of $T$ and $T'$ respectively.
Spoiler starts by selecting $v$. If Duplicator responds not with $v'$,
Spoiler applies the strategy of Lemma \ref{lem:new} and wins in
the next $d(T)$ moves. If Duplicator responds with $v'$,
Spoiler applies the strategy of Lemma \ref{lem:bothri*} and wins in
the next $\rk T+7$ moves. In any case Spoiler wins in
$1+\max\{d(T),\rk T+7\}=2\rk T+9$ moves.
\end{proof}

\begin{lemma}\label{lem:conndistdiam}
Let $T$ be a ranked tree and $G$ be either a tree of different diameter
or a connected non-tree. Then $D_0(T,G)\le2\rk T+10$.
\end{lemma}

\begin{proof}
If $G$ is a tree, then $d(T)+2$ moves are enough for Spoiler to win.
In this case, he selects a path of length $\min\{d(T),d(G)\}+1$
in the graph of larger diameter.

Suppose that $G$ is not a connected non-tree.
If $G$ has a cycle on at most $d(T)+2$ vertices, Spoiler selects it and wins.
Otherwise $G$ must have a cycle on at least $d(T)+3$ vertices.
Spoiler wins by selecting a path on $d(T)+2$ vertices of this cycle.
\end{proof}

\begin{lemma}\label{lem:nonranked}
Let $T$ be a ranked tree and $G$ be a non-ranked tree.
If $d(T)=d(G)$, then $D_0(T,G)\le2\rk T+9$.
\end{lemma}

\begin{proof}
Let $v$ and $c$ denote the central vertices of $T$ and $G$ respectively.
The tree in which Spoiler plays will be specified below. In the first
move Spoiler selects the central vertex of this tree. If Duplicator
responds not with the central vertex of the other tree, he loses in the next
$d(T)$ moves by Lemma \ref{lem:new}. Assume that she responds with the
central vertex.
Further play depends on which of three categories $G$ belongs to.
Let $k=\rk T$.
For any $w\in V(G)$ at distance $k$ from $c$, we will call
$G_c(w)$ an {\em apex\/} of $G_c$.

\case 1{$G_c$ has an apex $G_c(w)$ which is not a rooted subtree
of any of the four rooted trees in $R^*_0$}
Spoiler plays in $G$.
In the next $k$ moves he selects the path from $c$ to $w$.
Duplicator is enforced to select the path from $v$ to a vertex $u$
such that $T_v(u)\in R^*_0$. Spoiler is now able to win by selecting
at most 8 vertices of $G_c(w)$. The total number of moves does not
exceed $1+k+8=k+9$.

\case 2{$G$ has a vertex $w$ such that $B(G_c(w))$ properly
contains $B(H_w)$ for some $H_w\in R^*_i$, where $i=k-d(c,w)$}
Spoiler plays in $G$.
In the next $d(c,w)$ moves he selects the path from $c$ to $w$.
Let $u$ denote the vertex selected by Duplicator in response to $w$
and $F_u=T_v(u)$. Clearly, Duplicator must ensure the equality
$d(v,u)=d(c,w)$ and hence $F_u\in R^*_i$.

If $F_u$ and $H_w$ are not
isomorphic, then Spoiler restricts further play to $H_w$
following a continuous strategy. Of course, Duplicator
is enforced to play in $F_u$. Spoiler is able to win in the
next $i+7$ moves according to Lemma~\ref{lem:bothri*}.

Suppose now that $F_u$ and $H_w$ are isomorphic.
In the next move Spoiler selects a child of $w$ which is not
in $H_w$. Duplicator must respond with a child of $u$ in $F_u$.
Denote it by $x$ and let $y$ be the vertex of $H_w$ corresponding
to $x$ under the isomorphism from $F_u$ to $H_w$.
Recall that, by Lemma \ref{lem:divvsasym}, diverging trees are
asymmetric and therefore such isomorphism is unique.
In the next move Spoiler selects $y$. Duplicator must respond with $z$,
another child of $u$ in $F_u$. Note that $F_u(z)$ and $H_w(y)$
are not isomorphic since the latter is isomorphic to $F_u(x)$
but the former is not. From now on Spoiler restricts play to
$F_u(z)$ and $H_w(y)$ using the strategy of Lemma \ref{lem:bothri*},
and wins in the next $i+6$ moves. The total number of moves is at most
$1+d(c,w)+i+8=k+9$.

\case 3{Neither 1 nor 2}
Spoiler plays all the time in $T$. We will denote the vertices
selected by him in the next $k$ moves by $x_1,\ldots,x_k$ subsequently.
Let $y_1,\ldots,y_k$ denote the corresponding vertices selected in $G$
by Duplicator. Put also $x_0=v$ and $y_0=c$.
Spoiler will play so that $x_0,x_1,\ldots,x_k$
will be a path. Let $1\le i\le k$. Suppose that the preceding
$x_0,\ldots,x_{i-1}$ are already selected.
Assume that $T_v(x_{i-1})$ and $G_c(y_{i-1})$ are non-isomorphic
(note that this is so for $i=1$).
As we are not in Case 2, $x_{i-1}$ has a child $x$ such that
$T_v(x)\notin B(G_c(y_{i-1}))$. Spoiler takes this $x$ for $x_i$
thereby ensuring that $T_v(x_{i})$ and $G_c(y_{i})$ are non-isomorphic
again, whatever $y_i$ is selected by Duplicator.
The final stage of the game goes on non-isomorphic
$T_v(x_{k})$ and $G_c(y_{k})$. Spoiler selects all vertices of $T_v(x_{k})$.

Note that $T_v(x_{k})\in R^*_0$ and $G_c(y_{k})$ is an apex of $G$.
As we are not in Case 1, $G_c(y_{k})$ is a rooted subtree of some
$T^*_j\in R^*_0$. If $T^*_j=T_v(x_{k})$, $G_c(y_{k})$ must be a proper
subtree of $T_v(x_{k})$ and hence Spoiler has won. Otherwise, note that
$T_v(x_{k})$ cannot be a rooted subtree of $G_c(y_{k})$
by Property (Z4). Again, this is Spoiler's win. The total number of moves
equals $1+k+7=k+8$.

In any of the three cases Spoiler wins in $\max\{1+d(T),k+9\}=2k+9$ moves.
\end{proof}

Note that, if $T$ is a ranked tree of rank $k$, then Lemmas
\ref{lem:bothri}--\ref{lem:nonranked} provide Spoiler
with a winning strategy in the 0-alternation $\game_{2k+10}(T,G)$
whenever $G$ is a connected graph non-isomorphic to~$T$.

\begin{lemma}\label{lem:disconn}
Let $T$ be a ranked tree and $H$ be a disconnected graph.
Then $D_0(T,H)\allowbreak\le2\rk T+10$.
\end{lemma}

\begin{proof}
We distinguish two cases.

\case 1{No component of $H$ is isomorphic to $T$}

\subcase{1.1}{$H$ has a component $G$ such that Spoiler is able
to win $\game_{2k+10}(T,\allowbreak G)$ playing all the time in $G$}
Spoiler plays exactly this game.

\subcase{1.2}{$H$ has no such component}
In the first move Spoiler selects the central vertex of $T$.
Suppose that Duplicator's response is in a component $G$ of $H$.
By Lemmas \ref{lem:bothri}--\ref{lem:nonranked}, we are either
in the situation of Lemma \ref{lem:conndistdiam} (with $G$ a tree of
diameter $d(G)<d(T)$)
or in the situation of Lemma \ref{lem:nonranked} (namely, in Case 3).
In both situations Spoiler has a continuous winning strategy
for $\game_{2k+10}(T,G)$ allowing him to play all the time in $T$
starting from the central vertex. Spoiler applies it and wins as
Duplicator is enforced to stay in~$G$.

\case 2{$H$ has a component $T'$ isomorphic to $T$}
Spoiler plays in $H$. His first move is outside $T'$.
Let $x\in V(T)$ be Duplicator's response. Let $x'$ be the counterpart
of $x$ in $T'$ (recall that ranked trees are asymmetric and hence
$x'$ is determined uniquely). Denote the central vertices of $T$ and $T'$
by $v$ and $v'$ respectively. In the second move Spoiler
selects $v'$. If Duplicator responds not with $v$, Spoiler applies
the strategy of Lemma \ref{lem:new} and wins in the next $d(T)$ moves.
Assume that Duplicator responds with $v$. Starting from the third move,
Spoiler selects the vertices on the path between $v'$ and $x'$, one by one,
starting from a child of $v'$. If Duplicator follows the
path from $v$ to $x$, she loses as $x$ is already selected.
Assume that Duplicator deviates at some point, selecting a vertex $y$
not on the path, and let $y'$ be the vertex on the path between $v'$ and $x'$
selected in this round by Spoiler. Note that the rooted subtrees $T_v(y)$ and
$T'_{v'}(y')$ are non-isomorphic. Spoiler therefore can apply the
continuous strategy of Lemma \ref{lem:bothri*} and win in the next
$i+7$ moves, where $i=k-d(v,y)$.
The total number of moves is at most
$1+\max\{1+d(T),1+d(x,y)+(i+7)\}=2k+10$.
\end{proof}

Lemma \ref{lem:disconn} completes our analysis:
If $T$ is a ranked tree of rank $k$
and $G$ is an arbitrary graph non-isomorphic to $T$, then
we have a winning strategy for Spoiler in the 0-alternation
$\game_{2k+10}(T,G)$. By Proposition \ref{prop:ddd}, we conclude that
$D_0(T)\le 2\rk T+10$.

To complete the proof of Theorem \ref{thm:alt0}, let $T_i$ be a tree
of rank $i$ and order $N_i$ as in Lemma \ref{lem:minorder}.
We have $q_0(N_i;\mbox{\it trees})\le D_0(T_i)\le 2i+10\le2\log^*N_i+O(1)$,
the latter inequality due to Lemma~\ref{lem:minorder}.

\section{The succinctness function over trees: Lower bound}\label{s:tr_low}

Complementing the upper bound given by Theorem \ref{thm:tr_up}
we now prove a nearly tight lower bound on $\sftree$.

\begin{theorem}\label{thm:tr_low}
$\sftree\ge\log^*n-\log^*\log^*n-O(1)$.
\end{theorem}

It will be helpful to work with rooted trees. The first order language
for this class of structures has a constant $R$ for the root and the
parent-child relation
$P(x,y)$.  Let $T_v$ and $T'_u$ be rooted trees and suppose that
$T_v\keq T'_u$. By Proposition \ref{prop:loggames}, $T_v$ and $T'_u$
satisfy the same sentences of
quantifier rank $k$.  Then $T\keq T'$ for the underlying trees.
Indeed, take any sentence in the language for trees
and replace the adjacency $x\sim y$ with $P(x,y)\vee P(y,x)$.
We get a sentence with the same truth value in the language of
rooted trees.

Let $g(k)$ be the number of $\equiv_k$-equivalence classes
of rooted trees. Similarly to Lemma \ref{lem:ehrvk}, we have
$g(k)\le T(k+2+\log^*k)+O(1)$. Set
$$
U(k) = \sum_{i=0}^{g(k)-1} (k\,g(k))^i.
$$

\begin{lemma}\label{lem:main}
Let $T_v$ be a finite rooted tree. Then, for any $k\ge1$,
there exists a finite rooted tree $T'_u$ with at most $U(k)$
vertices such that $T_v\keq T'_u$.
\end{lemma}

\begin{proofof}{Theorem \ref{thm:tr_low}}
Consider an arbitrary tree $T$ of order $n$ and let $k=D(T)$.
Rooting it at an arbitrary vertex $v$, consider a rooted tree $T_v$.
Let $T'_u$ be as in Lemma \ref{lem:main}. Thus, we have
$T\keq T'$ and $|T'|\le U(k)$. By the choice of $k$, $T$ and $T'$
must be isomorphic. We therefore have
$$
n\le U(k)<\of{kg(k)}^{g(k)}\le T(k+\log^*k+4)+O(1),
$$
which implies $k\ge\log^*n-\log^*\log^*n-O(1)$.
\end{proofof}

Lemma \ref{lem:main} follows from a series of lemmas.

\begin{lemma}\label{lem:1}
Let $T_v$ be a rooted tree and $w$ a non-root vertex of $T_v$.
Suppose that $T_w' \keq T_v(w)$.
Let $T_v'$ be the result of replacing $T_v(w)$ by $T_w'$.
Then $T_v \keq T'_v$.
\end{lemma}

\begin{proof}
Duplicator wins the \EF\/ game on $T_v$, $T'_v$ by playing it on
$T_v(w)$, $T_w'$ (since the root is a constant symbol she automatically
plays root for root) and the identical vertices elsewhere.
\end{proof}

\begin{lemma}\label{lem:2}
Let $T_v$ be a rooted tree with $w_1,\ldots,w_s$ the children of the
root $v$, and $\ga_1,\ldots,\ga_s$ the $k$-Ehrenfeucht values of the
trees $T_v(w_i)$. Then the $k$-Ehrenfeucht value of $T$ is determined by
the $\ga_i$'s.
\end{lemma}

\begin{proof}
If $T_v$ and $T'_u$ have the same $\ga_1,\ldots,\ga_s$ we reach $T'_u$
from $T_v$ in $s$ applications of Lemma~\ref{lem:1}.
\end{proof}

\begin{lemma}\label{lem:4}
Suppose, in the notation of Lemma \ref{lem:2}, that some value $\ga$ appears
as $\ga_i$ more than $k$ times.  Let $T_v^-$ be $T_v$ but with only $k$ of
those subtrees. Then $T_v\keq T_v^-$.
\end{lemma}

\begin{proof}
The game has only $k$ moves so Spoiler cannot go in more
than $k$ of these subtrees.
\end{proof}

\begin{lemma}\label{lem:5}
If $T_v$ is a representative of a given $\keq$-equivalence class
with minimum possible order, then each vertex of $T_v$
has at most $kg(k)$ children.
\end{lemma}

\begin{proof}
This easily follows from Lemmas \ref{lem:4} and \ref{lem:2}
by induction on the depth.
\end{proof}

\begin{lemma}\label{lem:6}
If $T_v$ is a representative of a given $\keq$-equivalence class
with minimum possible order, then it has depth at most $g(k)-1$.
\end{lemma}

\begin{proof}
Take a longest path from the root to a leaf. If it
has more than $g(k)$ vertices, it contains two vertices $w$ and $u$ such that
$u$ is a descendant of $w$ and $T_v(u)\keq T_v(w)$. Replacing $T_v(w)$
by $T_v(u)$, we obtain a smaller tree in the same $\keq$-class.
\end{proof}

Lemma \ref{lem:main} immediately follows from Lemmas \ref{lem:5}
and~\ref{lem:6}.

\section{The smoothed succinctness function}\label{s:ssf}

Let $q(n)=q(n;\mbox{\it all})$ denote the succinctness function
for the class of all graphs. Since there are only finitely many
pairwise inequivalent sentences of bounded quantifier rank,
$q(n)\to\infty$ as $n\to\infty$. We will show that $q(n)$ grows
very slowly and, in a sense, irregularly. We first summarize
an information given by Theorems \ref{thm:hugegap} and~\ref{thm:tr_up}.

\begin{corollary}\label{cor:sf}
\mbox{}

\begin{enumerate}
\item
There is no general recursive function $f$ such that
$f(q(n))\ge n$ for all $n$.
\item
There is no general recursive function $l(n)$ such that $l(n)$ is
monotone nondecreasing, $l(n)\to\infty$ as $n\to\infty$, and
$l(n)\le q(n)$ for all $n$.
\item
$q(n)<\log^*n+5$.
\end{enumerate}
\end{corollary}

\begin{proof}
1)
Note $q(n)\le s(n)\le s_3(n)$. Now, if there were
a general recursive function $f$ such that $f(q(n))\ge n$, then we
would have $\max_{i\le s_3(n)}f(i)\ge n$ contradictory to
Theorem~\ref{thm:hugegap}.

2)
Assume that a such $l(n)$ exists. Let $f(m)$ be the first value of $i$
such that $l(i)>m$. Then $f(q(n))>n$ contradictory to Item~1.

3)
As any upper bound on $\sfclass$ is stronger if it
is proved for a smaller class of graphs, this item is
an immediate consequence of Theorem~\ref{thm:tr_up}.
\end{proof}

\begin{definition}
We define the {\em smoothed succinctness function\/} $q^*(n)$
(for quantifier rank) to be the least monotone nondecreasing
integer function bounding $q(n)$ from above, that is,
$q^*(n)=\max_{m\le n} q(m)$.
\end{definition}

\begin{theorem}\label{thm:ssf}
$\log^*n-\log^*\log^*n-O(1)<q^*(n)<\log^*n+5$.
\end{theorem}

\begin{proof}
Since the upper bound on $q(n)$ given by Corollary \ref{cor:sf} (3)
is monotone, this is a bound on $q^*(n)$ as well.
The lower bound is derivable from Lemma \ref{lem:ehrvk}.
This lemma states that $|\EHRV k|\le T(k+2+\log^*k)+c$ for a constant $c$.
Given $n>c+T(3)$, let $k$ be such that
$T(k+2+\log^*k)+c<n\le T(k+3+\log^*(k+1))+c$.
Assuming that $n$ is sufficiently large, we have
$k>\log^*n-\log^*\log^*n-4$.
According to Proposition \ref{prop:ddd}, at most $|\EHRV k|$
graphs are definable with quantifier rank at most $k$.
By the pigeonhole principle, there will be some
$m\leq |\EHRV k|+1\le n$ for which no graph of order precisely $m$ is defined
with quantifier rank at most $k$. We conclude that $q^*(n)\ge q(m)>k$
and hence $q^*(n)\ge \log^*n-\log^*\log^*n-2$.
\end{proof}

We defined $q^*(n)$ to be the ``closest'' to $q(n)$
monotone function. Notice that $q(n)$ itself lacks the monotonicity.

\begin{corollary}\label{cor:nonmon}
$q(i+1)<q(i)$ for infinitely many $i$.
\end{corollary}

\begin{proof}
Set $l(n)=\log^*n-\log^*\log^*n-2$.
We have just shown that $q^*(n)\ge l(n)$ for all $n$ large enough.
By Corollary \ref{cor:sf} (2), we have $q(n)<l(n)$ for infinitely many
$n$. For each such $n$, let $m_n<n$ be such that $q(m_n)\ge l(n)$.
Thus, $q(m_n)>q(n)$ and a desired $i$ must exist between $m_n$ and~$n$.
\end{proof}

For each non-negative integer $a$ and for $a=1/2$,
define $q_a(n)=\min_{|G|=n} D_a(G)$ and $q^*_a(n)=\max_{m\le n}q_a(m)$.
As easily seen, Corollary \ref{cor:sf} (1) holds true for $q_3(n)$ as well.
Note a strengthening of Corollary \ref{cor:sf} (3) that follows from a
result in another our paper. Let $G(n,p)$ denote a random graph
on $n$ vertices distributed so that each edge appears with probability
$p$ and all edges appear independently from each other.

\begin{theorem}\label{thm:rsa}
{\bf \cite{KPSV}}
With probability approaching 1 as $n$ goes to the infinity,
$$
D_3(G(n,n^{-1/4}))=\log^*n+O(1).
$$
\end{theorem}

\begin{corollary}
$q_3(n)\le\log^*n+O(1)$ and hence
$\log^*n-\log^*\log^*n-O(1)\le q_3^*(n)\le\log^*n+O(1)$.
\end{corollary}

\section{Depth vs.\ length}

\begin{theorem}\label{thm:dvsl}
$L(G)\le T(D(G)+\log^* D(G)+O(1))$.
\end{theorem}

\begin{proof}
Given an Ehrenfeucht value $\alpha$, let $l(\alpha)$ denote
the shortest length of a formula defining $\alpha$ in the sense of
Section \ref{s:game}.
Define $l(k)$ to be the maximum $l(\alpha)$ over
$\alpha\in\EHRV k$ and $l(k,s)$ the maximum $l(\alpha)$ over
$\alpha\in\EHRV{k,s}$. Of course, $l(k)=l(k,0)$.
As in Section \ref{s:game}, $f(k,s)=|\EHRV{k,s}|$.

It is not hard to see that $L(G)\le l(D(G))$ and therefore it suffices
to prove the bound $l(k)\le T(k+\log^*k+O(1))$ for all $k\ge2$.

On the account of Lemma \ref{lem:defehrv}, we have
$$
l(k,k)<18{k\choose 2}
$$
and
$$
l(k,s)\le f(k,s+1)(l(k,s+1)+10)
$$
if $s<k$. We will use these relations along with the bounds of
Lemma \ref{lem:ehrvks} for $f(k,s)$. Set $g(x)=x2^{x+1}$.
A simple inductive argument shows that
$$
f(k,s)\le 2^{g^{(k-s)}(9k^2)}\mbox{\ \ and\ \ }
l(k,s)\le g^{(k-s)}(9k^2).
$$
Since $g(x)\le 4^x$, we have
$l(k,0)\le T_4(k+2+\log^*k)\le T(k+\log^*k+O(1))$,
where $T_4$ stands for the variant of the tower function
built from 4's instead of 2's.
\end{proof}

\begin{remark}
Theorem \ref{thm:dvsl} generalizes to structures over
an arbitrary vocabulary. The proof requires only slight modifications.
\end{remark}

We now observe that the relationship between the optimum quantifier
rank and length of defining formulas is nearly tight.

\begin{theorem}\label{thm:tight}
There are infinitely many pairwise non-isomorphic graphs $G$ with
$L(G)\ge T(D(G)-6)-O(1)$.
\end{theorem}

\begin{proof}
The proof is given by a simple counting argument which can be naturally
presented in the framework of Kolmogorov complexity
(applications of Kolmogorov complexity for proving complexity-theoretic
lower bounds can be found in \cite{LVi}).

Denote the Kolmogorov complexity of a binary word $w$ by $K(w)$.
Let $\langle G\rangle$ denote the lexicographically first adjacency
matrix of a graph $G$. Define the Kolmogorov complexity of $G$
by $K(G)=K(\langle G\rangle)$. Notice that
$$
K(G)\le L(G)+O(1).
$$

By Theorem \ref{thm:tr_up}, there is a graph $G_n$ on $n$ vertices
with
\begin{equation}\label{eq:gnup}
D(G_n)<\log^*n+5.
\end{equation}
The bound $K(w)<k$ can hold for less than $2^k$ words. It follows that
for some $n\le 2^k$ we have $K(G)\ge k$ for all graphs $G$ on
$n$ vertices. For this particular $n$ we have
\begin{equation}\label{eq:gnlow}
L(G_n)\ge\log n-O(1).
\end{equation}
Combining \refeq{eq:gnup} and \refeq{eq:gnlow}, we see that
$G_n$ is as required.
\end{proof}

Of course, we could run the same argument directly with $L(G)$
in place of $K(G)$. An advantage of using the Kolmogorov complexity
is in avoiding estimation of the number of formulas of length at most~$k$.

In Section \ref{ss:pre} we showed that prenex formulas are sometimes
unexpectedly efficient in defining a graph. We are now able to show that,
nevertheless, they generally cannot be competetive against defining formulas
with no restriction on structure.
Let $\dpre G$ (resp.\ $\lpree G$) denote the minimum quantifier rank
(resp.\ length) of a closed prenex formula defining a graph~$G$.

\begin{theorem}
There are infinitely many pairwise non-isomorphic graphs $G$ with
$\dpre G\ge T(D(G)-8)$.
\end{theorem}

\begin{proof}
Let $G$ be as in Theorem \ref{thm:tight}. We have
$$
\lpree G\ge L(G)\ge T(D(G)-6)-O(1).
$$
On the other hand, by Lemma \ref{lem:lvsdprenex} we have
$$
\lpree G\le f(\dpre G)\mbox{,\ where\ }f(x)=O(x^24^{x^2}).
$$
It follows that
$$
\dpre G\ge\of{\frac1{\sqrt2}-o(1)}\sqrt{T(D(G)-7)}\ge T(D(G)-8),
$$
provided $D(G)$ (or the order of $G$) is sufficiently large.
\end{proof}

\section{Open questions}\label{s:open}
\mbox{}

\que
Let $D'(G)$ be the minimum quantifier rank of a first order sentence
distinguishing a graph $G$ from any non-isomorphic finite graph $G'$.
Clearly, $D'(G)\le D(G)$. Can the inequality be sometimes strict?

\que\label{que:alt}
Improve on the alternation number in Theorem \ref{thm:hugegap}.
The most interesting case is that of alternation number 0.
By the Ramsey theorem, Turing machines cannot be simulated
by 0-alternation formulas as this would contradict the unsolvability
of the halting problem. Thus, an intriguing question is how small
$s_0(n)$ and $q_0(n)$ can be.

\que
Classify the prefix classes with respect to solvability
of the finite satisfiability problem over graphs.
Such a classification does exist by the Gurevich classifiability
theorem \cite[section 2.3]{BGG}.
In particular, can the prefix
$\exists^*\forall^{O(1)}\exists^{O(1)}\forall^{O(1)}$
in Theorem \ref{thm:fgtundec}
be shortened to $\exists^*\forall^{O(1)}\exists^{O(1)}$?
Shortening to $\exists^*\forall^*$ is impossible due to
the Ramsey theorem.

Note that for digraphs the complete classification is known
(see \cite{BGG} and references there). In notation of Section \ref{ss:tfg},
the minimal undecidable classes for $\finsateq\calD$ are
$\forall^*\exists$,
$\forall\exists\forall^*$,
$\forall\exists\forall\exists^*$,
$\forall\exists^*\forall$,
$\exists^*\forall\exists\forall$,
$\exists^*\forall^{c+1}\exists$,
$\forall^{c+1}\exists^*$,
while the maximal decidable classes are
$\exists^*\forall^*$ and
$\exists^*\forall^c\exists^*$,
where $c=1$.
For $\finsat\calD$ the classification is the same but with $c=2$.
If we consider $\sateq\calD$ instead of $\finsateq\calD$ and
$\sat\calD$ instead of $\finsat\calD$, nothing in classification
changes. The reasons are that the maximal decidable classes have
the finite model property and that the undecidability of the minimal
undecidable classes is proved by reductions which preserve the finiteness
of models.

\que
How close to one another are $D_1(G)$ and $D_0(G)$?
At least, are they recursively linked?
The same question for $D(G)$ and $D_a(G)$
(any $a=o(n)$ is of interest).
How far apart from one another can be $D(G)$ and $D_1(G)$?

\que
Estimate the succinctness function $\sfclass$ for other classes of graphs
(in particular, graphs of bounded degree, planar graphs).

\que
Is $q(n)$ a non-recursive function? Is $D(G)$ an incomputable
function of graphs (T.~\L uczak)? Of course, the former implies the latter.
The same can be asked for $q_a(n)$ and $D_a(G)$ excepting $a\in\{0,1/2\}$
(see Theorem~\ref{thm:comput0}).

\que
We know that $q^*_3(n)=(1+o(1))\log^*n$. The cases of alternation number
0, 1, and 2 are open.

\que
$|q(n+1)-q(n)|=O(1)$? Note
$g(n+1)-g(n)\le 1$ but this difference is negative infinitely often
by Corollary \ref{cor:nonmon}.

\que
Can one construct a family of graphs $G$ as in Theorem \ref{thm:tight}
explicitly?

\subsection*{Acknowledgments}

We thank Nikolai Vereshchagin for a useful discussion.

\end{document}